\documentclass[a4paper,twoside,11pt]{article}
\usepackage{amsfonts}
\usepackage{amssymb}\usepackage{amsmath}
 \numberwithin{equation}{section}
\begin{document}

\title{The Cauchy problem of a periodic 2-component $\mu$-Hunter-Saxton system in Besov spaces}

\author{Jingjing Liu\footnote{E-mail: jingjing830306@163.com.}
\\Department of Mathematics and Information Science, \\ Zhengzhou
University of Light Industry , 450002 Zhengzhou, China}
\date{}
\maketitle

\noindent \textbf{Abstract}:
This paper is concerned with the local well-posedness and the
precise blow-up scenario for a periodic 2-component
$\mu$-Hunter-Saxton system in Besov spaces. Moreover, we state a new
global existence result to the system. Our obtained results for the
system improve considerably
earlier results.




\section{Introduction}
\par
Recently, a new 2-component system  was introduced by Zuo in
\cite{d-z} as follows:
\begin{equation}
\left\{\begin{array}{ll}
\mu(u)_{t}-u_{txx}=2\mu(u)u_{x}-2u_{x}u_{xx}-uu_{xxx}+\rho\rho_{x}&-\gamma_{1}u_{xxx},\\
&t > 0,\,x\in \mathbb{R},\\
\rho_{t}=(\rho u)_x+2\gamma_{2}\rho_{x}, &t > 0,\,x\in \mathbb{R},\\
u(0,x) = u_{0}(x),& x\in \mathbb{R}, \\
\rho(0,x) = \rho_{0}(x),&x\in \mathbb{R},\\
u(t,x+1)=u(t,x), & t \geq 0, x\in \mathbb{R},\\
\rho(t,x+1)=\rho(t,x), & t \geq 0, x\in \mathbb{R},\\ \end{array}\right. \\
\end{equation}
where $\mu(u)=\int_{\mathbb{S}}udx$ with
$\mathbb{S}=\mathbb{R}/\mathbb{Z}$ and $\gamma_{i}\in \mathbb{R},$
$i=1,2.$ By integrating both sides of the first equation in the
system (1.1) over the circle $\mathbb{S}=\mathbb{R}/\mathbb{Z}$ and
using the periodicity of $u$, one obtain
$$\mu(u_{t})=\mu(u)_{t}=0.$$
This yields the following periodic 2-component $\mu$-Hunter-Saxton
system:
\begin{equation}
\left\{\begin{array}{ll}
-u_{txx}=2\mu(u)u_{x}-2u_{x}u_{xx}-uu_{xxx}+\rho\rho_{x}&-\gamma_{1}u_{xxx},\\
&t > 0,\,x\in \mathbb{R},\\
\rho_{t}=(\rho u)_x+2\gamma_{2}\rho_{x}, &t > 0,\,x\in \mathbb{R},\\
u(0,x) = u_{0}(x),& x\in \mathbb{R}, \\
\rho(0,x) = \rho_{0}(x),&x\in \mathbb{R},\\
u(t,x+1)=u(t,x), & t \geq 0, x\in \mathbb{R},\\
\rho(t,x+1)=\rho(t,x), & t \geq 0, x\in \mathbb{R},\\ \end{array}\right. \\
\end{equation}
with $\gamma_{i}\in \mathbb{R},$ $i=1,2$. This system is a
2-component generalization of the generalized Hunter-Saxton equation
obtained in \cite{k-l-m}. The author \cite{d-z} shows that this
system is both a bihamiltonian Euler equation and a bivariational
equation. Moreover, the geometric background of the system (1.1) has
been comprehensively studied by Escher in \cite{escher1} recently.

Obviously, (1.1) is equivalent to (1.2) under the condition
$\mu(u_{t})=\mu(u)_{t}=0.$ In this paper, we will study the system
(1.2) under the assumption $\mu(u_{t})=\mu(u)_{t}=0$.

For $\rho\equiv 0$ and $\gamma_{i}=0, i=1,2$ and replacing $t$ by
$-t,$ the system (1.2) reduces to the generalized Hunter-Saxton
equation (named $\mu$-Hunter-Saxton equation or $\mu$-Camassa-Holm
equation) as follows:
\begin{equation}
-u_{txx}=-2\mu(u)u_{x}+2u_{x}u_{xx}+uu_{xxx},
\end{equation}
which is obtained and studied in \cite{k-l-m}. Moreover, the
periodic $\mu$-Hunter-Saxton equation and the periodic
$\mu$-Degasperis-Procesi equation have also been studied in
\cite{fu, l-m-t} recently. It is worthy to note that the
$\mu$-Hunter-Saxton equation has a very closed relation with the
periodic Hunter-Saxton and Camassa-Holm equations \cite{k-l-m}.

For $\rho\not\equiv0$ , $\gamma_{i}=0, i=1,2,$ $\mu(u)=0$ and
replacing $t$ by $-t,$ the system (1.2) becomes to a 2-component
periodic Hunter-Saxton system. Its peakon solutions and the Cauchy
problem have been analysed and discussed in \cite{C-I} and
\cite{MW}, respectively.

The system (1.2) has been studied in \cite{l-y1} in Sobolev spaces
$H^{s}(\mathbb{S})\times H^{s-1}(\mathbb{S}),$ $s\geq 2$ recently.
The authors established the local well-posedness in
$H^{s}(\mathbb{S})\times H^{s-1}(\mathbb{S}),$ $s\geq 2,$ by Kato's
semigroup theory, derived the precise blow-up scenario, presented
some blow-up results for strong solutions and gave a global
existence result to the system. Inspired by the study of the local
well-posedness and blow-up criterion to the Camassa-Holm equation in
\cite{d,d1} and the local well-posedness and blow-up criterion to
the two-component Camassa-Holm equation in \cite{gui1, gui2}, we
will discuss the system (1.2) in Besov spaces. Our obtained local
well-posedness, blow-up criterion and global existence results for
the system improve considerably earlier results in \cite{l-y1}.
Moreover, a very interesting result in this paper is Lemma 3.3.
Using this lemma, we will give a explicit proof of the continuity of
solution with respect to the initial data when establish the local
well-posedness of the system (1.2) in Besov spaces. In my opinion,
this proof is new and necessary.

Our paper is organized as follows. In Section 2, we recall some
basic facts of periodic Besov spaces and the transport equation
theory. In Section 3, we establish the local well-posedness of the
initial value problem associated with the system (1.2). In Section
4, we derive the precise blow-up scenario of strong solution to the
system (1.2) and present a
new global existence result for strong solutions to the system (1.2) with certain initial profiles.\\

 \textbf{Notation} Given a Banach space $Z$, we denote its norm by
 $\|\cdot\|_{Z}$. Since all space of functions are over
 $\mathbb{S}$, for simplicity, we drop $\mathbb{S}$ in our notations
  if there is no ambiguity. Let $u^{(k)}$ stand for $k$th derivate of
  $u$ and let $\ast$ denote the convolution.

\section{Preliminaries}

In this section, we will recall some basic facts on periodic Besov
spaces and the transport equation theory. We refer to
\cite{a-b,d,d1,st} for the elementary properties of them. Here, we
only display some facts which will be used later.\\
\newline
\textbf{Proposition 2.1} (\cite{a-b,d,d1,st,g}). \ (Littlewood-Paley
decomposition) Let $B=\{\xi\in\mathbb{R}, \ |\xi|\leq\frac{4}{3}\}$
and $C=\{\xi\in\mathbb{R}, \frac{3}{4}\leq|\xi|\leq\frac{8}{3}\}.$
There exist two radial functions $\chi\in C_{c}^{\infty}(B)$ and
$\varphi \in C_{c}^{\infty}(C)$ such that
$$\chi(\xi)+\sum\limits_{q\geq 0}\varphi(2^{-q}\xi)=1, \ \ \forall
\ \ \xi\in \mathbb{R}.$$ For $u\in \mathcal {D}^{'}(\mathbb{S}),$ let
$$\Delta_{q}u=0 \ \text{for} \ q\leq-2, \ \ \ \
\Delta_{-1}u=\sum\limits_{\beta\in\mathbb{Z}}\chi(\beta)\widehat{u}_{\beta}e^{2\pi
i \beta x},$$
$$\Delta_{q}u=\sum\limits_{\beta\in\mathbb{Z}}\varphi(2^{-q}\beta)\widehat{u}_{\beta}e^{2\pi
i \beta x} \ \ \ \ \text{for} \ \ \  q\geq0$$ and
$$S_{q}u=\sum\limits_{-1\leq p\leq q-1}\Delta_{p}u.$$ A direct computation implies, for any $u\in \mathcal {D}^{'}(\mathbb{S})$ and $v\in \mathcal {D}^{'}(\mathbb{S})$, the following properties hold:
$$\Delta_{p}\Delta_{q}u\equiv 0\ \ \ \ \text{if} \ |p-q|\geq 2,$$
$$\Delta_{q}(S_{p-1}u\Delta_{p}v)\equiv0 \ \ \ \ \text{if} \ |p-q|\geq 5.$$ Moreover, $$\|\Delta_{q}u\|_{L^{p}}\leq C\|u\|_{L^{p}}$$ for some constant $C$ independent of $q.$\\
\newline
\textbf{Definition 2.1} (\cite{a-b,d,d1,g,st}). \ (Besov spaces) Let
$s\in \mathbb{R},$ $1\leq p,r\leq\infty.$ The periodic Besov spaces
$B_{p,r}^{s}(\mathbb{S})$ is defined by
$$B_{p,r}^{s}(\mathbb{S})=\{u\in \mathcal {D}^{\prime}(\mathbb{S});
\|u\|_{B_{p,r}^{s}(\mathbb{S})}<\infty\},$$ where
$$
 \|u\|_{B_{p,r}^{s}(\mathbb{S})}= \left\{\begin{array}{ll}
\left(\sum\limits_{q\in\mathbb{Z}}2^{qsr}\|\Delta_{q}u\|_{L^{p}}^{r}\right)^{\frac{1}{r}}, \ &for \ r<\infty,\\
 \sup\limits_{q\in\mathbb{Z}}2^{qs}\|\Delta_{q}u\|_{L^{p}}, &for \ r=\infty.\end{array}\right.\\$$
If $s=\infty,$ $B_{p,r}^{\infty}=\cap_{s\in\mathbb{R}}B_{p,r}^{s}.$
The Sobolev spaces correspond to $H^{s}=B_{2,2}^{s}.$\\
\newline
\textbf{Proposition 2.2}
(\cite{a-b,d,d1,g,st}). \ The following properties hold:\\
(1) Density: for $1\leq p,r\leq\infty,$ we have $\mathcal
{D}(\mathbb{S})\subset B_{p,r}^{s} \subset \mathcal
{D}^{\prime}(\mathbb{S}).$ Moreover, if $p,r < \infty$, then the set
of all trigonometric polynomial is dense in
$B_{p,r}^{s}(\mathbb{S}).$ \\
(2) Sobolev embedding: if $p_{1}\leq p_{2}$ and $r_{1}\leq r_{2}$,
then $B_{p_{1},r_{1}}^{s}\hookrightarrow
B_{p_{2},r_{2}}^{s-(\frac{1}{p_{1}}-\frac{1}{p_{2}})}.$ If
$s_{1}<s_{2},$ $1\leq p \leq+\infty$ and $1\leq r_{1},
r_{2}\leq+\infty,$ then
$B_{p,r_{2}}^{s_{2}}\hookrightarrow B_{p,r_{1}}^{s_{1}}.$  Moreover, if $r_{1}=r_{2}$, then the embedding is compact.\\
(3) Algebraic properties: $(B_{p,r}^{s} \ \text{is an algebra})
\Leftrightarrow (B_{p,r}^{s}\hookrightarrow
L^{\infty})\Leftrightarrow (s>\frac{1}{p} \ \text{or} \
(s\geq\frac{1}{p}, \ \text{r=1})).$\\
(4) Fatou property: if $(u^{n})_{n \in \mathbb{N}}$ is a bounded
sequence of $B_{p,r}^{s}$ which tends to $u$ in $\mathcal
{D}^{\prime}(\mathbb{S}),$ then $u\in B_{p,r}^{s}$ and
$$\|u\|_{B_{p,r}^{s}}\leq
\liminf\limits_{n\rightarrow \infty}\|u^{n}\|_{B_{p,r}^{s}}.$$ (5)
Complex interpolation: if $u\in B_{p,r}^{s}\cap
B_{p,r}^{\widetilde{s}}$ and $\theta \in [0,1],$ $1\leq
p,r\leq\infty,$ then $u\in B_{p,r}^{\theta
s+(1-\theta)\widetilde{s}},$ and $\|u\|_{B_{p,r}^{\theta
s+(1-\theta)\widetilde{s}}}\leq
\|u\|_{B_{p,r}^{s}}^{\theta}\|u\|_{B_{p,r}^{\widetilde{s}}}^{1-\theta}.$\\
(6) The lifting property: let $u\in
\mathcal{D}^{\prime}(\mathbb{S})$ and $\alpha\in\mathbb{R}.$ Then
$u\in B_{p,r}^{s}$ if and only if
$$\sum\limits_{\beta\neq 0}e^{2 \pi i \beta
x}(i\beta)^{\alpha}\widehat{u}_{\beta}\in B_{p,r}^{s-\alpha}.$$ (7)
Let $s>0.$ Then $u\in B_{p,r}^{s+1}$ if and only if $u$ is
differentiable a.e. and $u^{\prime}\in B_{p,r}^{s}.$\\
\newline
\textbf{Lemma 2.1} (\cite{d, d1, g}). \ Suppose that $(p,r)\in
[1,+\infty]^{2}$ and $s>-\min\{\frac{1}{p}, 1-\frac{1}{p}\}.$ Let
$v$ be a vectorfield such that $\partial_{x}v$ belongs to
$L^{1}([0,T];B_{p,r}^{s-1})$ if $s>1+\frac{1}{p}$ or to
$L^{1}([0,T];B_{p,r}^{\frac{1}{p}}\cap L^{\infty})$ otherwise.
Suppose also that $f_{0}\in B_{p,r}^{s},$ $g\in
L^{1}([0,T];B_{p,r}^{s})$ and that $f\in
L^{\infty}([0,T];B_{p,r}^{s})\cap C([0,T];\mathcal
{D}^{\prime}(\mathbb{S}))$ solves the following linear transport
equation
\begin{equation}
\left\{\begin{array}{l}
\partial_{t}f+v\partial_{x}f=g,\\
 f|_{t=0}=f_{0},\\
 f(t,x+1)=f(t,x).\end{array}\right.\tag{T}  \label{T}
\end{equation}
Then there exists a constant $C$ depending only on $s,p$,
such that the following statements hold for all $t\in [0,T]$:\\
(i).
\begin{equation*}
\|f(t)\|_{B^s_{p,r}}\leq \|f_0\|_{B^s_{p,r}}\,+\, \int_0^t
\|g(\tau)\|_{B^s_{p,r}}d\tau\,+\, C\int_0^t
V^{'}(\tau)\|f(\tau)\|_{B^s_{p,r}} d\tau,
\end{equation*}
or hence
$$
\|f\|_{B_{p,r}^{s}}\leq
e^{CV(t)}\left(\|f_{0}\|_{B_{p,r}^{s}}+\int_{0}^{t}e^{-CV(\tau)}\|g(\tau)\|_{B_{p,r}^{s}}d\tau\right),
$$
where
$$
V(t)=\left\{\begin{array}{ll} \int_{0}^{t}\|\partial_{x}v(\tau,
\cdot)\|_{B_{p,r}^{\frac{1}{p}}\cap
L^{\infty}}d \tau, \ \ \ \ \text{if} \ \ s<1+\frac{1}{p},\\
\int_{0}^{t}\|\partial_{x}v(\tau,\cdot)\|_{B_{p,r}^{s-1}}d\tau,
\ \ \ \ \ \ \ \ \ \text{if} \ \ s>1+\frac{1}{p} \ \text{or} \ \{s=1+\frac{1}{p} \ and \ r=1\}.\end{array}\right.\\
$$
(ii). If $f=c_{1}(v+c_{2})$ with $c_{1}, c_{2}\in\mathbb{R},$ then
for all $s>0$, the estimates in (i) hold with
$V(t)=\int_{0}^{t}\|\partial_{x}v(\tau)\|_{L^{\infty}}d\tau.$\\
(iii). If $r<+\infty,$ then $f\in C([0,T];B_{p,r}^{s}).$ If
$r=+\infty,$ then $f\in C([0,T];B_{p,1}^{s^{\prime}})$ for all
$s^{\prime}<s.$\\
\newline
\textbf{Lemma 2.2} (\cite{d,gui2}). \  Let $0<\sigma<1$. Suppose
that $f_{0}\in H^{\sigma},$ $g\in L^{1}([0,T];H^{\sigma}),$ $v,
\partial_{x}v\in L^{1}([0,T]; L^{\infty})$ and that $f\in
L^{\infty}([0,T];H^{\sigma})\cap C([0,T]; \mathcal
{D}^{\prime}(\mathbb{S}))$ solves the 1-dimensional linear transport
equation
\begin{equation*}
\left\{
\begin{array}{l}
\partial_{t}f+v\partial_{x}f=g,\\
 f|_{t=0}=f_{0},\\
 f(t,x+1)=f(t,x).
 \end{array}
 \right. \tag{T}  \label{T}
\end{equation*}
Then $f\in C([0,T];H^{\sigma}).$ More precisely, there exists a
constant $C$ depending only on $\sigma$ and such that the following
statement holds:
$$
\|f(t)\|_{H^{\sigma}}\leq\|f_{0}\|_{H^{\sigma}}+C\int_{0}^{t}
\|g(\tau)\|_{H^{\sigma}}d\tau+C\int_{0}^{t}\|f(\tau)\|_{H^{\sigma}}V^{\prime}(\tau)d\tau
$$
or hence $$\|f\|_{H^{\sigma}}\leq
e^{CV(t)}\left(\|f_{0}\|_{H^{\sigma}}+C\int_{0}^{t}\|g(\tau)\|_{H^{\sigma}}d\tau\right)$$
with
$V(t)=\int_{0}^{t}(\|v(\tau)\|_{L^{\infty}}+\|\partial_{x}v(\tau)\|_{L^{\infty}})d\tau.$\\
\newline
\textbf{Lemma 2.3} (\cite{d,d1}). \ Let $(p,p_{1},r)\in
[1,+\infty]^{3}$. Assume that
$s>-\min\{\frac{1}{p_{1}},\frac{1}{p^{\prime}}\}$ or
$s>-1-\min\{\frac{1}{p_{1}},\frac{1}{p^{\prime}}\}$ if
$\partial_{x}v=0$ with $p^{\prime}=(1-\frac{1}{p})^{-1}$. Let
$f_{0}\in B_{p,r}^{s}$ and $g\in L^{1} ([0,T];B_{p,r}^{s}).$  Let
$v$ be a time dependent vector field such that $v\in
L^{\rho}([0,T];B_{\infty, \infty}^{-M})$ for some $\rho>1,$ $M>0$
and $\partial_{x}v\in
L^{1}([0,T];B_{p_{1},\infty}^{\frac{1}{p_{1}}}\cap L^{\infty})$ if
$s<1+\frac{1}{p_{1}},$  and $\partial_{x}v\in
L^{1}([0,T];B_{p_{1},r}^{s-1})$ if $s>1+\frac{1}{p_{1}}$ or
$s=1+\frac{1}{p_{1}}$ and $r=1.$ Then the transport equations $(T)$
has a unique solution $f\in L^{\infty}([0,T]; B_{p,r}^{s})\cap
(\cap_{s^{\prime}<s}C([0,T];B_{p,1}^{s^{\prime}}))$ and the
inequalities of Lemma 2.1 hold. If, moreover, $r < \infty,$ then we
have $f\in C([0,T]; B_{p,r}^{s}).$\\
\newline
\textbf{Lemma2.4} (\cite{ch,d,g,gui2}). \ (1-D Moser-type estimates)
Assume that $1\leq
p,r\leq +\infty,$ the following estimates hold:\\
(1) for $s>0$, $$\|fg\|_{B_{p,r}^{s}}\leq
C(\|f\|_{B_{p,r}^{s}}\|g\|_{L^{\infty}}+\|g\|_{B_{p,r}^{s}}\|f\|_{L^{\infty}});$$
(2) for $s>0$, $$\|f\partial_{x}g\|_{H^{s}}\leq
C(\|f\|_{H^{s+1}}\|g\|_{L^{\infty}}+\|f\|_{L^{\infty}}\|\partial_{x}g\|_{H^{s}});$$
(3) for $s_{1}\leq \frac{1}{p},
s_{2}>\frac{1}{p}(s_{2}\geq\frac{1}{p} \ \text{if} \ r=1)$ and
$s_{1}+s_{2}>0$, $$\|fg\|_{B_{p,r}^{s_{1}}}\leq
C\|f\|_{B_{p,r}^{s_{1}}}\|g\|_{B_{p,r}^{s_{2}}},$$ where $C$ is
constant independent of $f$ and $g$.

\section{Local well-posedness}

In this section, we will establish the local well-posedness for the
Cauchy problem of the system (1.2) in Besov spaces and then get our
main result in $H^{s}\times H^{s-1}$, $s>\frac{3}{2}$.

Note that $\mu(u)_{t}=\mu(u_{t})=0.$ Then we let
$$\mu_{0}=\mu(u_{0})=\mu(u)=\int_{\mathbb{S}}u(t,x)dx.$$
We now provide the framework in which we shall reformulate the
system (1.2). We rewrite the system (1.2) as follows:
\begin{equation}
\left\{\begin{array}{ll}
u_{t}-(u+\gamma_{1})u_{x}=\partial_{x}(\mu-\partial_{x}^{2})^{-1}&(2\mu_{0}
u+\frac{1}{2}u_{x}^{2}+\frac{1}{2}\rho^{2}),\\
&t > 0,\,x\in \mathbb{R},\\
\rho_{t}-(u+2\gamma_{2})\rho_{x}=u_{x}\rho, &t > 0,\,x\in \mathbb{R},\\
u(0,x) = u_{0}(x),&x\in \mathbb{R}, \\
\rho(0,x) = \rho_{0}(x),&x\in \mathbb{R},\\
u(t,x+1)=u(t,x), & t \geq 0, x\in \mathbb{R},\\
\rho(t,x+1)=\rho(t,x), & t \geq 0, x\in \mathbb{R}.\\
\end{array}\right. \\
\end{equation}

If we denote $P(D)$ as the Fourier integral operator with the
Fourier multiplier $-\frac{2\pi i
\beta}{\delta(\beta)+4\pi^{2}\beta^{2}}$ with
$$\delta(\beta)=\left\{\begin{array}{ll}1, \ \ \ \ \beta=0,\\
0, \ \ \ \ \beta\neq0, \end{array}\right. \\$$ then the system (3.1)
equivalent to
\begin{equation}
\left\{\begin{array}{ll} u_{t}-(u+\gamma_{1})u_{x}=P(D)(2\mu_{0}
u+\frac{1}{2}u_{x}^{2}+\frac{1}{2}\rho^{2}),
&t > 0,\,x\in \mathbb{R},\\
\rho_{t}-(u+2\gamma_{2})\rho_{x}=u_{x}\rho, &t > 0,\,x\in \mathbb{R},\\
u(0,x) = u_{0}(x),&x\in \mathbb{R}, \\
\rho(0,x) = \rho_{0}(x),&x\in \mathbb{R},\\
u(t,x+1)=u(t,x), & t \geq 0, x\in \mathbb{R},\\
\rho(t,x+1)=\rho(t,x), & t \geq 0, x\in \mathbb{R}.\\
\end{array}\right. \\
\end{equation} Moreover, combining Proposition 2.2 (6) and
\begin{eqnarray*}
P(D)u=\sum\limits_{\beta\in\mathbb{Z}}e^{2\pi i \beta
x}\widehat{P(D)u}_{\beta}&=&\sum\limits_{\beta\in\mathbb{Z}}e^{2\pi
i \beta x}\cdot\left(-\frac{2\pi i
\beta}{\delta(\beta)+4\pi^{2}\beta^{2}}\right)\widehat{u}_{\beta}\\
&=&\frac{1}{2\pi}\sum\limits_{\beta\neq0}e^{2\pi i \beta
x}(i\beta)^{-1}\widehat{u}_{\beta},
\end{eqnarray*} we have if $u\in B_{p,r}^{s},$ then
$\|P(D)u\|_{B_{p,r}^{s+1}}\leq C\|u\|_{B_{p,r}^{s}}.$

On the other hand, integrating both sides of the first equation in
(1.2) with respect to $x$, we obtain
$$u_{tx}=-2\mu_{0}u+\frac{1}{2}u_{x}^{2}+uu_{xx}-\frac{1}{2}\rho^{2}+\gamma_{1}u_{xx}+a(t),$$
where
$$a(t)=2\mu(u)^{2}+\frac{1}{2}\int_{\mathbb{S}}(u_{x}^{2}+\rho^{2})dx.$$
Using the system (1.2), we have
$$\frac{d}{dt}\int_{\mathbb{S}}(u_{x}^{2}+\rho^{2})dx=0.$$ By $\mu(u)_{t}=\mu(u_{t})=0,$ we have
$$\frac{d}{dt}a(t)=0.$$
For convenience, we let
$$\mu_{1}:=\left(\int_{\mathbb{S}}(u_{x}^{2}+\rho^{2})dx\right)^{\frac{1}{2}}
=\left(\int_{\mathbb{S}}(u_{0,x}^{2}+\rho_{0}^{2})dx\right)^{\frac{1}{2}}$$
and write $a:=a(0)$ henceforth. Thus,
\begin{equation}
u_{tx}=-2\mu_{0}u+\frac{1}{2}u_{x}^{2}+uu_{xx}-\frac{1}{2}\rho^{2}+\gamma_{1}u_{xx}+a
\end{equation}
is a valid reformulation of the first equation in (1.2).\\
\newline
\textbf{Definition 3.1.} For $T>0,$ $s\in\mathbb{R}$ and $1\leq p,
r\leq+\infty$, we set\\
$E_{p,r}^{s}(T)=C([0,T];B_{p,r}^{s})\cap C^{1}([0,T];B_{p,r}^{s-1})
\ if \ r<+\infty,$\\
$E_{p,\infty}^{s}(T)=L^{\infty}([0,T];B_{p,\infty}^{s})\cap
Lip([0,T];B_{p,\infty}^{s-1})$ and
$E_{p,r}^{s}=\cap_{T>0}E_{p,r}^{s}(T).$

The local well-posedness result of the system (1.2) in $B_{p,r}^{s}$
and $E_{p,r}^{s}(T)$ can be stated as follows:\\
\newline
\textbf{Theorem 3.1.} Suppose that $1\leq p,r\leq+\infty$ and
$s>\max\{1+\frac{1}{p}, 2-\frac{1}{p}, \frac{3}{2}\}$ with $s\neq
2+\frac{1}{p}.$ Given $z_{0}=(u_{0}, \rho_{0})\in B_{p,r}^{s}\times
B_{p,r}^{s-1},$ then there exists a time $T>0$ and a unique solution
$z=(u,\rho)$ to the system (3.2) such that $z=z(\cdot,z_{0})\in
E_{p,r}^{s}(T)\times E_{p,r}^{s-1}(T).$ Moreover, the mapping
$z_{0}\rightarrow z:$ is continuous from $B_{p,r}^{s}\times
B_{p,r}^{s-1}$ into $C([0,T]; B_{p,r}^{s^{\prime}}\times
B_{p,r}^{s^{\prime}-1})\times C^{1}([0,T];
B_{p,r}^{s^{\prime}-1}\times B_{p,r}^{s^{\prime}-2})$ for every
$s^{\prime}<s$ when $r=+\infty$ and $s^{\prime}=s$ otherwise.

Uniqueness with respect to the initial data is an immediate
consequence of the following result.\\
\newline
\textbf{Lemma 3.1.} Let $1\leq p, r\leq +\infty$ and
$s>\max\{1+\frac{1}{p}, 2-\frac{1}{p}, \frac{3}{2}\}.$ Suppose that
$(u^{i}, \rho^{i})\in L^{\infty}([0,T]; B_{p,r}^{s}\times
B_{p,r}^{s-1})\cap C([0,T]; B_{p,r}^{s-1}\times B_{p,r}^{s-2}),$
$i=1,2,$ are two solutions of the system (3.2) with initial data
$(u_{0}^{i}, \rho_{0}^{i})\in B_{p,r}^{s}\times B_{p,r}^{s-1},$
$i=1,2,$ then for every $t\in [0,T]$ with
$\mu_{0}^{i}=\int_{\mathbb{S}}u_{0}^{i}(x)dx,$ $i=1,2:$
\newline
(1) if $r=1,$ $s>\max\{1+\frac{1}{p}, 2-\frac{1}{p}, \frac{3}{2}\}$
or $r\neq 1,$ $s>\max\{1+\frac{1}{p}, 2-\frac{1}{p}, \frac{3}{2}\}$
but $s\neq 2+\frac{1}{p}, 3+\frac{1}{p}$, then
\begin{align*}
&\|u^{1}-u^{2}\|_{B_{p,r}^{s-1}}+\|\rho^{1}-\rho^{2}\|_{B_{p,r}^{s-2}}\\
\leq \ &
e^{C\int_{0}^{t}(\|u^{1}\|_{B_{p,r}^{s}}+\|u^{2}\|_{B_{p,r}^{s}}
+\|\rho^{1}\|_{B_{p,r}^{s-1}}+\|\rho^{2}\|_{B_{p,r}^{s-1}}+|\mu_{0}^{1}|)d\tau}
\\
&\left(\|u_{0}^{1}-u_{0}^{2}\|_{B_{p,r}^{s-1}}+\|\rho_{0}^{1}-\rho_{0}^{2}\|_{B_{p,r}^{s-2}}
+C|\mu_{0}^{1}-\mu_{0}^{2}|\int_{0}^{t}\|u^{2}(\tau)\|_{B_{p,r}^{s}}d\tau\right),
\end{align*}
(2) if $r\neq 1$ and $s=2+\frac{1}{p}$, then
\begin{align*}
&\|u^{1}-u^{2}\|_{B_{p,r}^{s-1}}+\|\rho^{1}-\rho^{2}\|_{B_{p,r}^{s-2}}\\
\leq \ &C
e^{C\theta\int_{0}^{t}(\|u^{1}\|_{B_{p,r}^{s}}+\|u^{2}\|_{B_{p,r}^{s}}
+\|\rho^{1}\|_{B_{p,r}^{s-1}}+\|\rho^{2}\|_{B_{p,r}^{s-1}}+|\mu_{0}^{1}|)d\tau}
\\
&
\left(\|u^{12}_0\|_{B^{s-1}_{p,r}}+\|\rho^{12}_0\|_{B^{s-2}_{p,r}}+C|\mu_{0}^{1}-\mu_{0}^{2}|\int_{0}^{t}\|u^{2}(\tau)\|_{B_{p,r}^{s}}d\tau\right)^{\theta}
\times\\
&
(\|u^{1}(t)\|_{B^{s}_{p,r}}+\|u^{2}(t)\|_{B^{s}_{p,r}})^{1-\theta}\\
+&e^{C\int_{0}^{t}(\|u^{1}\|_{B_{p,r}^{s}}+\|u^{2}\|_{B_{p,r}^{s}}
+\|\rho^{1}\|_{B_{p,r}^{s-1}}+\|\rho^{2}\|_{B_{p,r}^{s-1}}+|\mu_{0}^{1}|)d\tau}
\\
&\left(\|u_{0}^{1}-u_{0}^{2}\|_{B_{p,r}^{s-1}}+\|\rho_{0}^{1}-\rho_{0}^{2}\|_{B_{p,r}^{s-2}}
+C|\mu_{0}^{1}-\mu_{0}^{2}|\int_{0}^{t}\|u^{2}(\tau)\|_{B_{p,r}^{s}}d\tau\right),
\end{align*}
(3) if $r\neq 1$ and $s= 3+\frac{1}{p}$, then
\begin{align*}
&\|u^{1}-u^{2}\|_{B_{p,r}^{s-1}}+\|\rho^{1}-\rho^{2}\|_{B_{p,r}^{s-2}}\\
\leq \
&e^{C\int_{0}^{t}(\|u^{1}\|_{B_{p,r}^{s}}+\|u^{2}\|_{B_{p,r}^{s}}
+\|\rho^{1}\|_{B_{p,r}^{s-1}}+\|\rho^{2}\|_{B_{p,r}^{s-1}}+|\mu_{0}^{1}|)d\tau}
\\
&\left(\|u_{0}^{1}-u_{0}^{2}\|_{B_{p,r}^{s-1}}+\|\rho_{0}^{1}-\rho_{0}^{2}\|_{B_{p,r}^{s-2}}
+C|\mu_{0}^{1}-\mu_{0}^{2}|\int_{0}^{t}\|u^{2}(\tau)\|_{B_{p,r}^{s}}d\tau\right)\\
+&C
e^{C\theta\int_{0}^{t}(\|u^{1}\|_{B_{p,r}^{s}}+\|u^{2}\|_{B_{p,r}^{s}}
+\|\rho^{1}\|_{B_{p,r}^{s-1}}+\|\rho^{2}\|_{B_{p,r}^{s-1}}+|\mu_{0}^{1}|)d\tau}
\\
&
\left(\|u^{12}_0\|_{B^{s-1}_{p,r}}+\|\rho^{12}_0\|_{B^{s-2}_{p,r}}
+C|\mu_{0}^{1}-\mu_{0}^{2}|\int_{0}^{t}\|u^{2}(\tau)\|_{B_{p,r}^{s}}d\tau\right)^{\theta}
\times\\
&
(\|\rho^{1}(t)\|_{B^{s-1}_{p,r}}+\|\rho^{2}(t)\|_{B^{s-1}_{p,r}})^{1-\theta}.
\end{align*}
\newline
\textbf{Proof}
Denote $u^{12}=u^{1}-u^{2}$, $\rho^{12}=\rho^{1}-\rho^{2}.$ It is
obvious that $$u^{12}\in L^{\infty}([0,T]; B_{p,r}^{s})\cap C
([0,T]; B_{p,r}^{s-1}), \ \ \rho^{12}\in L^{\infty}([0,T];
B_{p,r}^{s-1})\cap C ([0,T]; B_{p,r}^{s-2}),$$ and $(u^{12},
\rho^{12})$ solves the transport equation:
 \begin{align*}
\left\{\begin{array}{ll}\partial_{t}u^{12}-(u^{1}+\gamma_{1})\partial_{x}u^{12}=u^{12}\partial_{x}u^{2}+F(t,x)
&t > 0,\,x\in \mathbb{R},\\
\partial_{t}\rho^{12}-(u^{1}+2\gamma_{2})\partial_{x}\rho^{12}=u^{12}\partial_{x}\rho^{2}
+\rho^{12}\partial_{x}u^{1}+\rho^{2}\partial_{x}u^{12}, &t > 0,\,x\in \mathbb{R},\\
u^{12}(0,x) = u^{1}_{0}(x)-u^{2}_{0}(x)=u_{0}^{12}(x),&x\in \mathbb{R}, \\
\rho^{12}(0,x) = \rho^{1}_{0}(x)-\rho^{2}_{0}(x)=\rho_{0}^{12}(x),&x\in \mathbb{R},\\
u^{12}(t,x+1)=u^{12}(t,x), & t \geq 0, x\in \mathbb{R},\\
\rho^{12}(t,x+1)=\rho^{12}(t,x), & t \geq 0, x\in \mathbb{R},\\
\end{array}\right.
\end{align*}
where
$$F(t,x)=P(D)(2\mu_{0}^{1}u^{12}+2(\mu_{0}^{1}-\mu_{0}^{2})u^{2}
+\frac{1}{2}\partial_{x}u^{12}\partial_{x}(u^{1}+u^{2})+\frac{1}{2}\rho^{12}(\rho^{1}+\rho^{2})).$$

(1) If $r=1,$ $s>\max\{1+\frac{1}{p}, 2-\frac{1}{p}, \frac{3}{2}\}$
or $r\neq 1,$ $s>\max\{1+\frac{1}{p}, 2-\frac{1}{p}, \frac{3}{2}\}$
but $s\neq 2+\frac{1}{p}, 3+\frac{1}{p},$ noting that for $w\in
B_{p,r}^{s}$ with $s>1+\frac{1}{p},$ then
$\|\partial_{x}w\|_{B_{p,r}^{{\frac{1}{p}}}\cap{L^{\infty}}}\leq
C\|w\|_{B_{p,r}^{s}}.$ Applying Lemma 2.1 and the fact that
$\|\partial_{x}w\|_{B_{p,r}^{s-3}}\leq
C\|\partial_{x}w\|_{B_{p,r}^{s-2}}\leq C\|w\|_{B_{p,r}^{s}},$ we
have
\begin{align}
& \nonumber
e^{-C\int_{0}^{t}\|u^{1}(\tau)\|_{B_{p,r}^{s}}d\tau}\|u^{12}(t)\|_{B_{p,r}^{s-1}}\leq
\|u_{0}^{12}\|_{B_{p,r}^{s-1}}\\
 \ &\ \ \ \ \ \ \
+\int_{0}^{t}e^{-C\int_{0}^{\tau}\|u^{1}(\tau^{\prime})\|_{B_{p,r}
^{s}}d\tau^{\prime}}\cdot
\left(\|u^{12}\partial_{x}u^{2}\|_{B_{p,r}^{s-1}}+\|F(\tau)\|_{B_{p,r}^{s-1}}\right)d\tau
\end{align}
and
\begin{align}
& \nonumber
e^{-C\int_{0}^{t}\|u^{1}(\tau)\|_{B_{p,r}^{s}}d\tau}\|\rho^{12}(t)\|_{B_{p,r}^{s-2}}\\
\nonumber  \leq \ &
\|\rho_{0}^{12}\|_{B_{p,r}^{s-2}}+\int_{0}^{t}e^{-C\int_{0}^{\tau}\|u^{1}(\tau^{\prime})\|_{B_{p,r}
^{s}}d\tau^{\prime}}\\
& \cdot
\left(\|u^{12}\partial_{x}\rho^{2}\|_{B_{p,r}^{s-2}}+\|\rho^{12}\partial_{x}u^{1}\|_{B_{p,r}^{s-2}}+
\|\rho^{2}\partial_{x}u^{12}\|_{B_{p,r}^{s-2}}\right)d\tau.
\end{align}
For $s>1+\frac{1}{p},$ $B_{p,r}^{s-1}$ is an algebra according to
Proposition 2.3 (3), so we have
$$\|u^{12}\partial_{x}u^{2}\|_{B_{p,r}^{s-1}}\leq \|u^{12}\|_{B_{p,r}^{s-1}}\|\partial_{x}u^{2}\|_{B_{p,r}^{s-1}}\leq
\|u^{12}\|_{B_{p,r}^{s-1}}\|u^{2}\|_{B_{p,r}^{s}}.$$ By the property
of $P(D),$ we have
\begin{align*}
&\|P(D)(2\mu_{0}^{1}u^{12}+2(\mu_{0}^{1}-\mu_{0}^{2})u^{2})\|_{B_{p,r}^{s-1}}\\
\leq \ &C\|2\mu_{0}^{1}u^{12}+2(\mu_{0}^{1}-\mu_{0}^{2})u^{2}\|_{B_{p,r}^{s-2}}\\
\leq \ &
C(|\mu_{0}^{1}|\|u^{12}\|_{B_{p,r}^{s-1}}+|\mu_{0}^{1}-\mu_{0}^{2}|\|u^{2}\|_{B_{p,r}^{s}}).
\end{align*}
Moreover, note that $B_{p,r}^{s-2}$ is an algebra with
$s-2>\frac{1}{p}.$ If $s-2\leq \frac{1}{p},$ then combining
$s>\max\{1+\frac{1}{p}, \frac{3}{2}\}$ and Lemma 2.4, we get
\begin{align*}
\|P(D)(\partial_{x}u^{12}\partial_{x}(u^{1}
+u^{2}))\|_{B_{p,r}^{s-1}}\leq & \
\|\partial_{x}u^{12}\partial_{x}(u^{1}
+u^{2})\|_{B_{p,r}^{s-2}} \\
\leq& \
C\|\partial_{x}u^{12}\|_{B_{p,r}^{s-2}}\left(\|\partial_{x}u^{1}\|_{B_{p,r}^{s-1}}
+\|\partial_{x}u^{2}\|_{B_{p,r}^{s-1}}\right)\\
\leq& \
C\|u^{12}\|_{B_{p,r}^{s-1}}\left(\|u^{1}\|_{B_{p,r}^{s}}+\|u^{2}\|_{B_{p,r}^{s}}\right),
\end{align*}
\begin{align*}
\|P(D)(\rho^{12}(\rho^{1}+\rho^{2}))\|_{B_{p,r}^{s-1}}\leq& \
C\|\rho^{12}\|_{B_{p,r}^{s-2}}\left(\|\rho^{1}\|_{B_{p,r}^{s-1}}+\|\rho^{2}\|_{B_{p,r}^{s-1}}\right).
\end{align*}
The inequalities above imply:
\begin{align}
\nonumber
\|u^{12}\partial_{x}u^{2}\|_{B_{p,r}^{s-1}}&+\|F(\tau)\|_{B_{p,r}^{s-1}}
\leq
C(\|u^{12}\|_{B_{p,r}^{s-1}}+\|\rho^{12}\|_{B_{p,r}^{s-2}})\cdot\\
\nonumber &  \ \ \ (\|u^{1}\|_{B_{p,r}^{s}}+\|u^{2}\|_{B_{p,r}^{s}}+
\|\rho^{1}\|_{B_{p,r}^{s-1}}+\|\rho^{2}\|_{B_{p,r}^{s-1}}+|\mu_{0}^{1}|)\\
 & \ \ \ \ \ \ \ \ \ \ \ \ \ \ \ \ \ \ \ \ \ \ \ \ \ \ \ \
\ \ \ \ \ \ \ \ +C|\mu_{0}^{1}-\mu_{0}^{2}|\|u^{2}\|_{B_{p,r}^{s}}.
\end{align}
While thanks to Lemma 2.4, we have
$$\|u^{12}\partial_{x}\rho^{2}\|_{B_{p,r}^{s-2}}\leq C\|u^{12}\|_{B_{p,r}^{s-1}}\|\partial_{x}\rho^{2}\|_{B_{p,r}^{s-2}}
\leq C\|u^{12}\|_{B_{p,r}^{s-1}}\|\rho^{2}\|_{B_{p,r}^{s-1}},$$
$$\|\rho^{12}\partial_{x}u^{1}\|_{B_{p,r}^{s-2}}\leq C\|\rho^{12}\|_{B_{p,r}^{s-2}}\|\partial_{x}u^{1}\|_{B_{p,r}^{s-1}}
\leq C\|\rho^{12}\|_{B_{p,r}^{s-2}}\|u^{1}\|_{B_{p,r}^{s}},$$ and
$$\|\rho^{2}\partial_{x}u^{12}\|_{B_{p,r}^{s-2}}\leq C\|\rho^{2}\|_{B_{p,r}^{s-1}}\|\partial_{x}u^{12}\|_{B_{p,r}^{s-2}}
\leq C\|\rho^{2}\|_{B_{p,r}^{s-1}}\|u^{12}\|_{B_{p,r}^{s-1}}.$$ It
then follows that
\begin{align}
&\nonumber
\|u^{12}\partial_{x}\rho^{2}\|_{B_{p,r}^{s-2}}+\|\rho^{12}\partial_{x}u^{1}\|_{B_{p,r}^{s-2}}+
\|\rho^{2}\partial_{x}u^{12}\|_{B_{p,r}^{s-2}}\\
\leq \ &
C(\|u^{12}\|_{B_{p,r}^{s-1}}+\|\rho^{12}\|_{B_{p,r}^{s-2}})(\|u^{1}\|_{B_{p,r}^{s}}+\|\rho^{2}\|_{B_{p,r}^{s-1}}).
\end{align}
Thus, combining (3.4)-(3.7), we have
\begin{align*}
&e^{-C\int_{0}^{t}\|u^{1}(\tau)\|_{B_{p,r}^{s}}d\tau}\left(\|u^{12}(t)\|_{B_{p,r}^{s-1}}
+\|\rho^{12}(t)\|_{B_{p,r}^{s-2}}\right)\\
\leq \ & \|u_{0}^{12}\|_{B_{p,r}^{s-1}}
+\|\rho_{0}^{12}\|_{B_{p,r}^{s-2}}+C\int_{0}^{t}
|\mu_{0}^{1}-\mu_{0}^{2}|\|u^{2}(\tau)\|_{B_{p,r}^{s}}d\tau\\
&+C\int_{0}^{t}e^{-C\int_{0}^{\tau}\|u^{1}(\tau^{\prime})\|_{B_{p,r}^{s}}d\tau^{\prime}}
\left(\|u^{12}(\tau)\|_{B_{p,r}^{s-1}}
+\|\rho^{12}(\tau)\|_{B_{p,r}^{s-2}}\right)\\
&\cdot\left(\|u^{1}\|_{B_{p,r}^{s}}+\|u^{2}\|_{B_{p,r}^{s}}
+\|\rho^{1}\|_{B_{p,r}^{s-1}}+\|\rho^{2}\|_{B_{p,r}^{s-1}}+|\mu_{0}^{1}|\right)d\tau.
\end{align*}
That is
$$w(t)\leq v(t)+C\int_{0}^{t} w(\tau) u(\tau) d\tau$$
with
$$w(t)= e^{-C\int_{0}^{t}\|u^{1}(\tau)\|_{B_{p,r}^{s}}d\tau}\left(\|u^{12}(t)\|_{B_{p,r}^{s-1}}
+\|\rho^{12}(t)\|_{B_{p,r}^{s-2}}\right),$$
$$v(t)=\|u_{0}^{12}\|_{B_{p,r}^{s-1}}
+\|\rho_{0}^{12}\|_{B_{p,r}^{s-2}}+C|\mu_{0}^{1}-\mu_{0}^{2}|\int_{0}^{t}\|u^{2}(\tau)\|_{B_{p,r}^{s}}d\tau,$$
and
$$u(t)=\|u^{1}\|_{B_{p,r}^{s}}+\|u^{2}\|_{B_{p,r}^{s}}
+\|\rho^{1}\|_{B_{p,r}^{s-1}}+\|\rho^{2}\|_{B_{p,r}^{s-1}}+|\mu_{0}^{1}|.$$
This completes the proof of (1) by applying Gronwall's inequality.

(2) If $r\neq 1$ and $s= 2+\frac{1}{p}$, we will use the
interpolation method to deal with it. Indeed, if we choose
$s_1\in(\max\{1+\frac{1}{p}, 2-\frac{1}{p},\frac 3 2\}-1,s-1)$,
$s_2\in (s-1,s)$ and $\theta=\frac{s_2-(s-1)}{s_2-s_1} \in(0,1)$,
then $s-1=\theta s_1+(1-\theta )s_2$. According to Proposition 2.3
(5), we have
\begin{align*}
&\|u^{12}(t)\|_{B^{s-1}_{p,r}}\\
\leq \ & \|u^{12}(t)\|^{\theta}_{B^{s_1}_{p,r}}\|u^{12}(t)\|^{1-\theta}_{B^{s_2}_{p,r}}\\
\leq \ &
(\|u^{1}(t)\|_{B^{s_2}_{p,r}}+\|u^{2}(t)\|_{B^{s_2}_{p,r}})^{1-\theta}\|u^{12}(t)\|^{\theta}_{B^{s_1}_{p,r}}.
\end{align*}
Since $s_{1}+1>\max\{1+\frac{1}{p}, 2-\frac{1}{p},\frac 3 2\}$ and
$s_{1}+1<s=2+\frac{1}{p},$ the estimate in case (1) for
$||u^{12}(t)||_{B^{s_1}_{p,r}}$ holds. On the other hand, thanks to
$s-2=\frac{1}{p}<1+\frac{1}{p}$, we have (3.5) holds. Consequently,
the estimate for $||\rho^{12}(t)||_{B^{s-2}_{p,r}}$ in case (1) can
also hold true. Hence, we can get the desired result.

(3) For the critical case $s= 3+\frac{1}{p}$, noting that $s-1=
2+\frac{1}{p}>1+\frac{1}{p}$, we have (3.4) holds. So the estimate
for $\|u^{12}(t)\|_{B^{s-1}_{p,r}}$ in case (1) holds here. The left
proof is very similar to that of case (2). Therefore, we complete
our proof of Lemma 3.1.

Next, we will construct the approximate solutions to (3.2).\\
\newline
\textbf{Lemma 3.2.} Let $u_{0}, \rho_{0}, p, r$ and $s$ be as in the
statement of Theorem 3.1. Assume that $u^{0}=\rho^{0}=0.$ Then there
exists a unique sequence of smooth functions $(u^{n},
\rho^{n})_{n\in N}\in C(\mathbb{R}^{+}; B_{p,r}^{\infty}\times
B_{p,r}^{\infty})$ solving the following linear transport equation
by induction:
\[ \ \ \ \ \left\{\begin{array}{l}
\partial_{t}u^{n+1}-(u^{n}+\gamma_{1})\partial_{x}u^{n+1}=F^{n}(t,x),\\
\partial_{t}\rho^{n+1}-(u^{n}+2\gamma_{2})
\partial_{x}\rho^{n+1}=\rho^{n}\partial_{x}u^{n},\\
u^{n+1}(0,x) = u^{n+1}_{0}(x)=S_{n+1}u_{0}, \\
\rho^{n+1}(0,x) = \rho^{n+1}_{0}(x)=S_{n+1}\rho_{0}, \\
u^{n+1}(t,x+1)=u^{n+1}(t,x),\\
\rho^{n+1}(t,x+1)=\rho^{n+1}(t,x),\\
\end{array}\right. \tag{Tn}  \label{Tn} \]
where $F^{n}(t,x)=P(D)(2\mu_{0}^{n+1}
u^{n}+\frac{1}{2}(\partial_{x}u^{n})^{2}+\frac{1}{2}(\rho^{n})^{2})$
with $\mu_{0}^{n+1}=\int_{\mathbb{S}}u_{0}^{n+1}dx.$ Moreover, there
exists a positive $T$ such that the solutions satisfy: \\
(i) $(u^{n}, \rho^{n})_{n\in N}$ is uniformly bounded in
$E_{p,r}^{s}(T)\times E_{p,r}^{s-1}(T).$\\
(ii) $(u^{n}, \rho^{n})_{n\in N}$ is a Cauchy sequence in $C([0,T];
B_{p,r}^{s-1}\times B_{p,r}^{s-2}).$\\
\newline
\textbf{Proof}
For convenience, we assume that $r\neq 1$ here. In fact, Theorem 3.2
corresponds to $p=r=2.$ Since all the data $S_{n+1}u_{0}$ and
$S_{n+1}\rho_{0}$ belong to $B_{p,r}^{\infty}$, Lemma 2.3 enables us
to show by induction that for all $n\in N$, the equation $(T_{n})$
has a unique global solution which belongs to $C(\mathbb{R}^{+},
B_{p,r}^{\infty}\times B_{p,r}^{\infty})$. Note that
$$|\mu_{0}^{n+1}|\leq\int_{\mathbb{S}}|u_{0}^{n+1}|dx\leq \|S_{n+1}u_{0}\|_{L^{\infty}}\leq
 \|u_{0}\|_{L^{\infty}}.$$ Similar to the
proof of Lemma 3.1, by $s>\max\{1+\frac{1}{p}, 2-\frac{1}{p},
\frac{3}{2}\}$ with $s\neq 2+\frac{1}{p},$ we have the following
inequalities for all $n\in N$,
\begin{align}
& \nonumber
e^{-C\int_{0}^{t}\|u^{n}(\tau)\|_{B_{p,r}^{s}}d\tau}\|u^{n+1}(t)\|_{B_{p,r}^{s}}\\
\nonumber  \leq \
&\|u_{0}\|_{B_{p,r}^{s}}+\frac{C}{2}\int_{0}^{t}e^{-C\int_{0}^{\tau}\|u^{n}(\tau^{\prime})\|_{B_{p,r}
^{s}}d\tau^{\prime}}\\
& \ \ \ \ \ \ \ \ \ \ \ \ \ \ \  \ \ \ \  \ \ \ \ \cdot
\left(\|u^{n}\|_{B_{p,r}^{s}}+\|u^{n}\|_{B_{p,r}^{s}}^{2}+\|\rho^{n}\|_{B_{p,r}^{s-1}}^{2}\right)d\tau
\end{align}
and
\begin{align}
& \nonumber
e^{-C\int_{0}^{t}\|u^{n}(\tau)\|_{B_{p,r}^{s}}d\tau}\|\rho^{n+1}(t)\|_{B_{p,r}^{s-1}}\\
\leq \ &
\|\rho_{0}\|_{B_{p,r}^{s-1}}+\frac{C}{2}\int_{0}^{t}e^{-C\int_{0}^{\tau}\|u^{n}(\tau^{\prime})\|_{B_{p,r}
^{s}}d\tau^{\prime}}
\|u^{n}\|_{B_{p,r}^{s}}\|\rho^{n}\|_{B_{p,r}^{s-1}}d\tau.
\end{align}
Hence we have
\begin{align}
&\nonumber e^{-C\int_{0}^{t}\|u^{n}(\tau)\|_{B_{p,r}^{s}}d\tau}
\left(\|u^{n+1}(t)\|_{B_{p,r}^{s}}+\|\rho^{n+1}(t)\|_{B_{p,r}^{s-1}}\right)\\
\nonumber\leq \ &
\left(\|u_{0}\|_{B_{p,r}^{s}}+\|\rho_{0}\|_{B_{p,r}^{s-1}}\right)
+\frac{C}{2}\int_{0}^{t}e^{-C\int_{0}^{\tau}\|u^{n}(\tau^{\prime})\|_{B_{p,r}
^{s}}d\tau^{\prime}}\\
&  \ \ \ \ \ \ \ \ \ \ \ \ \ \ \ \cdot
\left(\|u^{n}\|_{B_{p,r}^{s}}+\|\rho^{n}\|_{B_{p,r}^{s-1}}\right)
\left(\|u^{n}\|_{B_{p,r}^{s}}+\|\rho^{n}\|_{B_{p,r}^{s-1}}+1\right)d\tau.
\end{align}
Denoting $l^{n}(t)=
||u^n(t)||_{B^s_{p,r}}+||\rho^n(t)||_{B^{s-1}_{p,r}}$,  $L=
||u_0||_{B^s_{p,r}}+||\rho_0||_{B^{s-1}_{p,r}},$ we have
$$l^{n+1}(t)\leq
e^{C\int_{0}^{t}\|u^{n}(\tau)\|_{B_{p,r}^{s}}d\tau}\left(L+\frac{C}{2}\int_{0}^{t}e^{-C\int_{0}^{\tau}\|u^{n}(\tau^{\prime})\|_{B_{p,r}
^{s}}d\tau^{\prime}}l^{n}(l^{n}+1)d\tau\right).$$ Let us choose a
$T>0$ such that $T<\min\{\frac{1}{4CL}, \frac{1}{2C}\}$ and prove by
induction that for all $t\in [0,T]$
\begin{equation}
l^{n}(t)\leq\frac{2L}{1-4CLt}.
\end{equation}
Note that $$\int_{\tau}^{t}\|u^{n}(\tau^{\prime})\|_{B_{p,r}
^{s}}d\tau^{\prime}\leq -\frac{1}{2C}\ln\frac{1-4CLt}{1-4CL\tau}$$
with $l^{n}(t)\leq\frac{2L}{1-4CLt}.$ A direct computation implies
$$l^{n+1}(t)\leq\frac{2L}{1-4CLt}.$$
Therefore, $(u^{n}, \rho^{n})_{n\in N}$ is uniformly bounded in
$C([0,T]; B_{p,r}^{s}\times B_{p,r}^{s-1}).$ Using the fact that
$B_{p,r}^{s-1}$ with $s>1+\frac{1}{p}$ is an algebra, together with
Lemma 2.4, one can see that
$$(u^{n}+\gamma_{1})\partial_{x}u^{n+1}, \ \
P(D)(2\mu_{0}^{n+1}u^{n}+\frac{1}{2}(\partial_{x}u^{n})^{2}+\frac{1}{2}(\rho^{n})^{2})$$
are uniformly bounded in $C([0,T]; B_{p,r}^{s-1}),$  and
$$(u^{n}+2\gamma_{2})\partial_{x}\rho^{n+1}, \ \ \rho^{n}\partial_{x}u^{n}$$
are uniformly bounded in $C([0,T]; B_{p,r}^{s-2}).$ Hence using the
equations $(T_{n}),$ we have $$(\partial_{t}u^{n+1},
\partial_{t}\rho^{n+1})\in C([0,T]; B_{p,r}^{s-1}\times
B_{p,r}^{s-2})$$ are uniformly bounded, which yields that the
sequence $(u^{n}, \rho^{n})_{n\in N}$ is uniformly bounded in
$E_{p,r}^{s}(T)\times E_{p,r}^{s-1}(T).$\\

Next, we show that $(u^{n}, \rho^{n})_{n\in N}$ is a Cauchy sequence
in $C([0,T]; B_{p,r}^{s-1}\times B_{p,r}^{s-2}).$ In fact, according
to the equations $(T_{n}),$ we obtain that, for all $m,n\in N$
\begin{align}
&\nonumber(\partial_{t}-(u^{n+m}+\gamma_{1})\partial_{x})(u^{n+m+1}-u^{n+1})=(u^{n+m}-u^{n})\partial_{x}u^{n+1}\\
&\nonumber+P(D)(2\mu_{0}^{n+m+1}(u^{n+m}-u^{n}))+P(D)(2(\mu_{0}^{n+m+1}-\mu_{0}^{n+1})u^{n})\\
&\nonumber+\frac{1}{2}P(D)((\partial_{x}u^{n+m}-\partial_{x}u^{n})(\partial_{x}u^{n+m}+\partial_{x}u^{n}))\\
& +\frac{1}{2}P(D)((\rho^{n+m}-\rho^{n})(\rho^{n+m}+\rho^{n}))
\end{align}
and
\begin{align}
&\nonumber(\partial_{t}-(u^{n+m}+2\gamma_{2})\partial_{x})(\rho^{n+m+1}-\rho^{n+1})=(u^{n+m}-u^{n})\partial_{x}\rho^{n+1}\\
&
+\rho^{n+m}\partial_{x}(u^{n+m}-u^{n})+(\rho^{n+m}-\rho^{n})\partial_{x}u^{n}.
\end{align}
Applying Lemma 2.1, together with the fact that $B_{p,r}^{s-1}$ is
an algebra and the property of the operator $P(D)$, we get for
$t\in[0,T]$
\begin{align}
&\nonumber
e^{-C\int_{0}^{t}\|u^{n+m}(\tau)\|_{B_{p,r}^{s}}d\tau}\|(u^{n+m+1}-u^{n+1})(t)\|_{B_{p,r}^{s-1}}\\
\nonumber  \leq \ &
\|u^{n+m+1}_{0}-u_{0}^{n+1}\|_{B_{p,r}^{s-1}}+C\int_{0}^{t}e^{-C\int_{0}^{\tau}\|u^{n+m}(\tau^{\prime})\|_{B_{p,r}
^{s}}d\tau^{\prime}}\cdot\\
&\nonumber
[\|u^{n+m}-u^{n}\|_{B_{p,r}^{s-1}}(\|u^{n+1}\|_{B_{p,r}^{s}}+1
+\|u^{n}(\tau)\|_{B_{p,r}^{s}}+\|u^{n+m}(\tau)\|_{B_{p,r}^{s}})\\
&
+\|\rho^{n+m}-\rho^{n}\|_{B_{p,r}^{s-2}}(\|\rho^{n}\|_{B_{p,r}^{s-1}}
+\|\rho^{n+m}\|_{B_{p,r}^{s-1}})+|\mu_{0}^{n+m+1}-\mu_{0}^{n+1}|\|u^{n}\|_{B_{p,r}^{s}}]d\tau
\end{align}
and
\begin{align}
&\nonumber
e^{-C\int_{0}^{t}\|u^{n+m}(\tau)\|_{B_{p,r}^{s}}d\tau}\|(\rho^{n+m+1}-\rho^{n+1})(t)\|_{B_{p,r}^{s-2}}\\
 \nonumber  \leq \ &
\|\rho^{n+m+1}_{0}-\rho_{0}^{n+1}\|_{B_{p,r}^{s-2}}+C\int_{0}^{t}e^{-C\int_{0}^{\tau}\|u^{n+m}(\tau^{\prime})\|_{B_{p,r}
^{s}}d\tau^{\prime}}\\
&\nonumber
[\|u^{n+m}-u^{n}\|_{B_{p,r}^{s-1}}(\|\rho^{n+1}\|_{B_{p,r}^{s-1}}
+\|\rho^{n+m}\|_{B_{p,r}^{s-1}})\\
&
+\|\rho^{n+m}-\rho^{n}\|_{B_{p,r}^{s-2}}\|u^{n}\|_{B_{p,r}^{s}}]d\tau.
\end{align}
By Proposition 2.1 and Definition 2.3, we have
$$u_{0}^{n+m+1}-u_{0}^{n+1}=S_{n+m+1}u_{0}-S_{n+1}u_{0}=\sum\limits_{q=n+1}^{n+m}\triangle_{q}u_{0},$$
$$\rho_{0}^{n+m+1}-\rho_{0}^{n+1}=S_{n+m+1}\rho_{0}-S_{n+1}\rho_{0}=\sum\limits_{q=n+1}^{n+m}\triangle_{q}\rho_{0},$$
Moreover,
\begin{eqnarray*}
||\sum\limits_{q=n+1}^{n+m} \Delta_q u_0||_{B^{s-1}_{p,r}} \leq
C2^{-n}||u_0||_{B^s_{p,r}},
\end{eqnarray*}
and
$$||\sum\limits_{q=n+1}^{n+m} \Delta_q \rho_0||_{B^{s-2}_{p,r}}\leq C2^{-n}||\rho_0||_{B^{s-1}_{p,r}},$$
see \cite{yan} for detailed computations. Since $(u^{n},
\rho^{n})_{n\in N}$ is uniformly bounded in $E_{p,r}^{s}(T)\times
E_{p,r}^{s-1}(T)$, combining (3.14)-(3.15), we get a constant
$C_{T}$ independent of $n,m$ such that for all $t\in[0,T]$
\begin{equation}
h_{n+1}^{m}(t)\leq
C_{T}\left(2^{-n}+\int_{0}^{t}h_{n}^{m}(\tau)d\tau+|\mu_{0}^{n+m+1}-\mu_{0}^{n+1}|\right)
\end{equation}
with $$h_{n}^{m}(t)=\|(u^{n+m}-u^{n})(t)\|_{B_{p,r}^{s-1}}
+\|(\rho^{n+m}-\rho^{n})(t)\|_{B_{p,r}^{s-2}}.$$ Arguing by
induction with respect to the index $n$, one can easily prove that
\begin{align*}
&h^m_{n+1}(t) \\
\leq \ & C_{T}\left(2^{-n}\sum\limits_{k=0}^n \frac{(2T
C_T)^k}{k!}+\frac{(T C_T)^{n+1}}{(n+1)!}
+\sum\limits_{k=0}^n|\mu_{0}^{m+n-k+1}-\mu_{0}^{n-k+1}|\frac{(C_{T}T)^{k}}{k!}\right)\\
\leq \ & \big(C_{T}\sum\limits_{k=0}^n \frac{(2T
C_T)^k}{k!}\big)2^{-n}+C_T \frac{(T
C_T)^{n+1}}{(n+1)!}\\
& \ \ \ \ \ \ \ \ \ \ \ \ \ \ \ \ \ \ \ \ \ \ \ \ \ \ \ \ \
+C_T\sum\limits_{k=0}^n|\mu_{0}^{m+n-k+1}-\mu_{0}^{n-k+1}|\frac{(C_{T}T)^{k}}{k!},
\end{align*}
which implies that $(u^{n}, \rho^{n})_{n\in N}$ is a Cauchy sequence
in $C([0,T]; B_{p,r}^{s-1}\times B_{p,r}^{s-2}).$ This completes the
proof of Lemma 3.2.

Following the proof of Theorem 4.24 in \cite{bcd}, we obtain the
following result, which is crucial in the proof of the continuity of
solution with respect to the initial data.\\
\newline
\textbf{Lemma 3.3.} Denote $\overline{\mathbb{N}}=\mathbb{N}\cup
\{\infty\}.$ Suppose that $(p,r)\in[1,\infty]^{2},$ $r<+\infty,$
$s>1+\frac{1}{p}$ or $s\geq 1+\frac{1}{p}, \ r=1.$ Given a sequence
$\{a^{n}\}_{n\in\overline{\mathbb{N}}}$ of periodic continuous
bounded functions on $[0,T]\times \mathbb{S}$ with
$\partial_{x}a^{n}\in C([0,T]; B_{p,r}^{s-1})$ and for some
$\alpha(t)\in L^{1}([0,T])$
$$\|\partial_{x}a^{n}\|_{B_{p,r}^{s-1}}\leq\alpha(t), \ \ \ \
\text{for all} \ \ t\in [0,T], \ \ n\in \overline{\mathbb{N}}.$$
Assume that $\{v^{n}\}_{n\in\overline{\mathbb{N}}}\in
L^{\infty}([0,T]; B_{p,r}^{s-1})$ is the solution of
\begin{equation}
\left\{\begin{array}{ll}
\partial_{t}v^{n}+a^{n}\partial_{x}v^{n}=f,\\
 v^{n}|_{t=0}=v_{0},\\
 v^{n}(t,x+1)=v^{n}(t,x),\end{array}\right.\\
\end{equation} with $v_{0}\in B_{p,r}^{s-1},$ $f\in C([0,T],
B_{p,r}^{s-1}).$ If $a^{n}\rightarrow a^{\infty}$ in
$L^{1}([0,T];B_{p,r}^{s-1})$ as $n\rightarrow \infty,$ then the
sequence $\{v^{n}\}_{n\in \mathbb{N}}$ tends to $v^{\infty}$ in
$C([0,T], B_{p,r}^{s-1})$ as $n\rightarrow \infty.$\\
\newline
\textbf{Proof}
We first consider $v_{0}\in B_{p,r}^{s}$ and $f\in C([0,T],
B_{p,r}^{s}).$ Note that $r<\infty.$ By Lemma 2.1 (iii), we have
$\{v^{n}\}_{n\in\overline{\mathbb{N}}}\in C([0,T]; B_{p,r}^{s})$ in
this particular case. From (3.17), we get
$$\partial_{t}(v^{n}-v^{\infty})+a^{n}\partial_{x}(v^{n}-v^{\infty})=(a^{\infty}-a^{n})\cdot
\partial_{x}v^{\infty}.$$ Applying Lemma 2.1, we have
\begin{align*}
\|v^{n}-v^{\infty}\|_{B_{p,r}^{s-1}}\leq &
\int_{0}^{t}e^{C\int_{\tau}^{t}\|\partial_{x}a^{n}(\tau^{\prime})\|_{B_{p,r}^{s-1}}d\tau^{\prime}}\|(a^{\infty}-a^{n})\cdot
\partial_{x}v^{\infty}\|_{B_{p,r}^{s-1}}d\tau\\
\leq
&\int_{0}^{t}e^{C\int_{\tau}^{t}\alpha(\tau^{\prime})d\tau^{\prime}}\|a^{\infty}-a^{n}\|_{B_{p,r}^{s-1}}\|
\partial_{x}v^{\infty}\|_{B_{p,r}^{s-1}}d\tau\\
\leq
&\int_{0}^{t}e^{C\int_{\tau}^{t}\alpha(\tau^{\prime})d\tau^{\prime}}\|a^{\infty}-a^{n}\|_{B_{p,r}^{s-1}}\|v^{\infty}\|_{B_{p,r}^{s}}d\tau.
\end{align*}
Since $a^{n}\rightarrow a^{\infty}$ in $L^{1}([0,T];B_{p,r}^{s-1})$
as $n\rightarrow \infty,$ we have $v^{n}\rightarrow v^{\infty}$ in
$C([0,T];B_{p,r}^{s-1})$ as $n\rightarrow \infty.$

Now, we will discuss the general case $v_{0}\in B_{p,r}^{s-1},$
$f\in C([0,T], B_{p,r}^{s-1}).$ For all $n\in \overline{\mathbb{N}}$
and $q\in \mathbb{N},$ we consider the following equation:
\begin{equation}
\left\{\begin{array}{ll}
\partial_{t}v_{q}^{n}+a^{n}\partial_{x}v_{q}^{n}=S_{q}f,\\
 v_{q}^{n}|_{t=0}=S_{q}v_{0},\\
 v_{q}^{n}(t,x+1)=v_{q}^{n}(t,x).\end{array}\right.\\
\end{equation}
On one hand, since all the date $S_{q}v_{0}$ and $S_{q}f$ belongs to
$B_{p,r}^{\infty},$ the step above implies
\begin{equation}
v_{q}^{n}\rightarrow v_{q}^{\infty} \ \ \text {in} \
C([0,T];B_{p,r}^{s-1})\ \ \text{as} \ n\rightarrow \infty.
\end{equation}
On the other hand, for $n\in \overline{\mathbb{N}}$ and $q\in
\mathbb{N},$  subtracting (3.18) from (3.17) gives
$$
\left\{\begin{array}{ll}
\partial_{t}(v^{n}-v_{q}^{n})+a^{n}\partial_{x}(v^{n}-v_{q}^{n})=f-S_{q}f,\\
 (v^{n}-v_{q}^{n})|_{t=0}=v_{0}-S_{q}v_{0},\\
 (v^{n}-v_{q}^{n})(t,x+1)=(v^{n}-v_{q}^{n})(t,x).\end{array}\right.\\
$$
It follows that
\begin{align*}
&\|v^{n}-v_{q}^{n}\|_{B_{p,r}^{s-1}}\\
\leq \ &
e^{C\int_{0}^{t}\alpha(\tau)d\tau}\cdot\left(\|v_{0}-S_{q}v_{0}\|_{B_{p,r}^{s-1}}+
\int_{0}^{t}e^{-C\int_{0}^{\tau}\alpha(\tau^{\prime})d\tau^{\prime}}\|f-S_{q}f\|_{B_{p,r}^{s-1}}d\tau\right).
\end{align*}
By the definition of $S_{q}$ in Proposition 2.1 and the Lebesgue
dominated convergence theorem, we have for all
$n\in\overline{\mathbb{N}},$
\begin{equation}
v_{q}^{n}\rightarrow v^{n} \ \ \text {in} \ C([0,T];B_{p,r}^{s-1})\
\ \text{as} \ q\rightarrow \infty.
\end{equation}
Note that
$$\|v^{n}-v^{\infty}\|_{B_{p,r}^{s-1}}\leq\|v^{n}-v_{q}^{n}\|_{B_{p,r}^{s-1}}+\|v_{q}^{n}-v_{q}^{\infty}\|_{B_{p,r}^{s-1}}
+\|v_{q}^{\infty}-v^{\infty}\|_{B_{p,r}^{s-1}}.$$ For fixed $q$
large enough, letting $n$ tend to infinity, then combining (3.19)
and (3.20), we have the desired result.

Next, we give the proof of Theorem 3.1.\\
\newline
\textbf{Proof}
According to Lemma 3.2 (ii), we have that $z^{n}=(u^{n},
\rho^{n})_{n\in N}$ converges to some function $z=(u,\rho)\in
C([0,T]; B_{p,r}^{s-1}\times B_{p,r}^{s-2}).$ Next, we will prove
that $z=(u,\rho)$ satisfies Theorem 3.1.

Firstly, we will claim that $z=(u,\rho)$ is indeed a solution of the
system (3.2). Obviously, $z=(u,\rho)$ satisfies (3.2) in the sence
of $\mathcal {D}^{\prime}([0,T]\times\mathbb{R}).$ Combining (i) and
(ii) in Lemma 3.2 and using the interpolation estimate (5) in
Proposition 2.3, we have that $z^{n}=(u^{n}, \rho^{n})_{n\in N}$ is
a Cauchy consequence in $C([0,T]; B_{p,r}^{s^{\prime}}\times
B_{p,r}^{s^{\prime}-1}),$ for any $s^{\prime}< s$. Moreover,
$$z^n \rightarrow z,\,\ as \,\ n\to \infty,\,\ in\,\ C([0,T];
 B^{s^{\prime}}_{p,r}\times B^{s^{\prime}-1}_{p,r}). $$ Therefore,
$$(u+\gamma_{1})u_{x}+P(D)(2\mu_{0}
u+\frac{1}{2}u_{x}^{2}+\frac{1}{2}\rho^{2}) \ \ \ \text{and} \ \ \ \
(u+2\gamma_{2})\rho_{x}+u_{x}\rho$$ is continuous to $z=(u,\rho)$ in
$C([0,T]; B_{p,r}^{s^{\prime}-1}\times B_{p,r}^{s^{\prime}-2}).$
Taking limit in $(T_n)$, we can see that $z$ solves the system (3.2)
in the sense of $C([0,T]; B^{s^{\prime}-1}_{p,r})\times C([0,T];
B^{s^{\prime}-2}_{p,r})$ for all $s^{\prime}<s$. Furthermore,
combining
$$\int_{\mathbb{S}}u^{n}(x)dx=\int_{\mathbb{S}}u_{0}^{n}(x)dx\rightarrow
 \int_{\mathbb{S}}u_{0}(x)dx=\mu_{0}$$ and
 $$\int_{\mathbb{S}}u^{n}(x)dx\rightarrow \int_{\mathbb{S}}u(x)dx$$
 as $n\rightarrow \infty,$ we know that $u$ satisfies $\mu(u)_{t}=0.$

Secondly, we will prove that $z\in E_{p,r}^{s}(T)\times
E_{p,r}^{s-1}(T).$ Lemma 3.2 and Proposition 2.2 (4) guarantee that
$z=(u, \rho)$ belongs to $L^{\infty}([0,T];B_{p,r}^{s}\times
B_{p,r}^{s-1}).$ It follows that the right-hand side of the equation
$$u_{t}-(u+\gamma_{1})u_{x}=P(D)(2\mu_{0}
u+\frac{1}{2}u_{x}^{2}+\frac{1}{2}\rho^{2})$$ belongs to
$L^{\infty}([0,T];B_{p,r}^{s})$ and the right-hand side of the
equation
$$\rho_{t}-(u+2\gamma_{2})\rho_{x}=u_{x}\rho$$
belongs to $L^{\infty}([0,T];B_{p,r}^{s-1}).$ By Lemma 2.1 (iii), we
have that $z=(u, \rho)\in C([0,T]; B_{p,r}^{s^{\prime}}\times
B_{p,r}^{s^{\prime}-1})$ for any $s^{\prime}\leq s.$ Using the
system (3.2) again, we have $z\in E_{p,r}^{s}(T)\times
E_{p,r}^{s-1}(T).$

Thirdly, we will prove the continuity of solution with respect to
the initial data. At first, the continuity with respect to the
initial data in $$ C([0,T]; B^{s^{\prime}}_{p,r}\times
B^{s^{\prime}-1}_{p,r})\cap C^1{([0,T]; B^{s^{\prime}-1}_{p,r}\times
B^{s^{\prime}-2}_{p,r})}, \ \ \forall\, s^{\prime}<s$$ can be
obtained by Lemma 3.1 and a simple interpolation argument. Then we
will prove that the continuity holds true up to index $s.$ By the
argument before, we know there is a $B_{p,r}^{s}\times
B_{p,r}^{s-1}$-neighborhood $B_{z_{0}}$ of $z_{0}=(u_{0},\rho_{0})$
and some $T>0$ such that for any $v_{0}\in B_{z_{0}},$ the system
(3.2) with initial data $v_{0}$ has a solution $v\in C([0,T];
B^{s}_{p,r}\times B^{s-1}_{p,r})\cap C^1{([0,T]; B^{s-1}_{p,r}\times
B^{s-2}_{p,r})}.$ For $n\in \overline{\mathbb{N}},$ consider a
sequence of data $z_{0}^{n}=(u_{0}^{n},\rho_{0}^{n})\in B_{z_{0}}$
satisfing $z_{0}^{n}\rightarrow z_{0}^{\infty}:=z_{0}$ in
$B^{s}_{p,r}\times B^{s-1}_{p,r}.$ Then we have the corresponding
solutions $z^{n}=(u^{n}, \rho^{n})\in C([0,T]; B^{s}_{p,r}\times
B^{s-1}_{p,r})\cap C^1{([0,T]; B^{s-1}_{p,r}\times B^{s-2}_{p,r})}$
satisfy
\[(E_n) \ \ \ \ \left\{\begin{array}{l}
u^{n}_{t}-(u^{n}+\gamma_{1})u^{n}_{x}=P(D)(2\mu_{0}^{n}
u^{n}+\frac{1}{2}(u^{n}_{x})^{2}+\frac{1}{2}(\rho^{n})^{2}),\\
\rho^{n}_{t}-(u^{n}+2\gamma_{2})
\rho^{n}_{x}=u^{n}_{x}\rho^{n},\\
u^{n}(0,x) = u^{n}_{0}(x), \\
\rho^{n}(0,x) = \rho^{n}_{0}(x), \\
u^{n}(t,x+1)=u^{n}(t,x),\\
\rho^{n}(t,x+1)=\rho^{n}(t,x).\\
\end{array}\right.\]
Next, we will prove
$$z^{n}=(u^{n},\rho^{n})\rightarrow z=(u,\rho)\ \  \
\text{in} \ \ C([0,T]; B^{s}_{p,r}\times B^{s-1}_{p,r}) \ \ \
\text{as} \ \ n\rightarrow\infty.$$ Since $u^{n}\rightarrow u \ \
\text{in} \ C([0,T]; B^{s-1}_{p,r}),$ it suffices to prove that
$$u^{n}_{x}\rightarrow u_{x} \ \ \text{in} \ \
C([0,T]; B^{s-1}_{p,r})\ \ \ \text{and} \ \ \ \rho^{n}\rightarrow
\rho \ \ \text{in} \ \ C([0,T]; B^{s-1}_{p,r}).$$ For this purpose,
differentiating the first equation in $(E_{n})$ with respect to $x$,
we have
\[\left\{\begin{array}{l}
(u^{n}_{x})_{t}-(u^{n}+\gamma_{1})(u^{n}_{x})_{x}=F^{n},\\
\rho^{n}_{t}-(u^{n}+2\gamma_{2})\rho^{n}_{x}=u^{n}_{x}\rho^{n},\\
u^{n}(0,x) = u^{n}_{0}(x), \\
\rho^{n}(0,x) = \rho^{n}_{0}(x), \\
u^{n}(t,x+1)=u^{n}(t,x),\\
\rho^{n}(t,x+1)=\rho^{n}(t,x),\\
\end{array}\right.\]
where $F^{n}=(u^{n}_{x})^{2}+\partial_{x}P(D)(2\mu_{0}^{n}
u^{n}+\frac{1}{2}(u^{n}_{x})^{2}+\frac{1}{2}(\rho^{n})^{2})$ with
$\mu_{0}^{n}=\int_{\mathbb{S}}u_{0}^{n}(x)dx.$ Taking
$u^{n}_{x}=w^{n}+v^{n}$ and $\rho^{n}=f^{n}+g^{n},$ we have
$$\left\{\begin{array}{l}
w^{n}_{t}-(u^{n}+\gamma_{1})w^{n}_{x}=F,\\
f^{n}_{t}-(u^{n}+2\gamma_{2})f^{n}_{x}=u_{x}\rho,\\
w^{n}(0,x) = u_{0x}, \\
f^{n}(0,x) = \rho_{0}, \\
w^{n}(t,x+1)=w^{n}(t,x),\\
f^{n}(t,x+1)=f^{n}(t,x),
\end{array}\right.
\text{and} \ \ \left\{\begin{array}{l}
v^{n}_{t}-(u^{n}+\gamma_{1})v^{n}_{x}=F^{n}-F,\\
g^{n}_{t}-(u^{n}+2\gamma_{2})g^{n}_{x}=u^{n}_{x}\rho^{n}-u_{x}\rho,\\
v^{n}(0,x) =u^{n}_{0x}-u_{0x}, \\
g^{n}(0,x) =\rho_{0}^{n}-\rho_{0}, \\
w^{n}(t,x+1)=w^{n}(t,x),\\
f^{n}(t,x+1)=f^{n}(t,x),
\end{array}\right.$$
where $F=u_{x}^{2}+\partial_{x}P(D)(2\mu_{0}
u+\frac{1}{2}u_{x}^{2}+\frac{1}{2}\rho^{2}).$ By the first system
above and Lemma 3.3, we have
$$w^{n}\rightarrow w=u_{x} \ \ \text{in} \ \
C([0,T]; B^{s-1}_{p,r})\ \ \ \text{and} \ \ \ f^{n}\rightarrow \rho
\ \ \text{in} \ \ C([0,T]; B^{s-1}_{p,r}).$$ By Lemma 3.2 and $z\in
E_{p,r}^{s}(T)\times E_{p,r}^{s-1}(T),$ we obtain that there is a
positive constant $M$ such that $\|z^{n}\|_{B^{s}_{p,r}\times
B^{s-1}_{p,r}}\leq M$ and $\|z\|_{B^{s}_{p,r}\times
B^{s-1}_{p,r}}\leq M.$ By the second system above and Lemma 2.1, we
have
\begin{equation}
\|v^{n}\|_{B^{s-1}_{p,r}}\leq
e^{CMT}\left(\|u^{n}_{0x}-u_{0x}\|_{B^{s-1}_{p,r}}+\int_{0}^{t}\|F^{n}-F\|_{B^{s-1}_{p,r}}d\tau\right),
\end{equation}
and
\begin{equation}
\|g^{n}\|_{B^{s-1}_{p,r}}\leq
e^{CMT}\left(\|\rho_{0}^{n}-\rho_{0}\|_{B^{s-1}_{p,r}}+\int_{0}^{t}\|u^{n}_{x}\rho^{n}-u_{x}\rho\|_{B^{s-1}_{p,r}}d\tau\right).
\end{equation}
Noticing that for $s>1+\frac{1}{p},$ $B_{p,r}^{s-1}$ is an algebra,
we have
$$\|(u^{n}_{x})^{2}-u_{x}^{2}\|_{B^{s-1}_{p,r}}\leq\|u^{n}_{x}+u_{x}\|_{B^{s-1}_{p,r}}
\|u^{n}_{x}-u_{x}\|_{B^{s-1}_{p,r}}\leq 2M
\|u^{n}_{x}-u_{x}\|_{B^{s-1}_{p,r}}.$$ By the property of $P(D),$ we
have
\begin{align*}
&\|\partial_{x}P(D)(2\mu_{0}^{n}u^{n}-2\mu_{0}
u)\|_{B^{s-1}_{p,r}}\\
\leq \ &\|2\mu_{0}^{n}u^{n}-2\mu_{0}
u\|_{B^{s-1}_{p,r}}\\
\leq \ &
2|\mu_{0}^{n}-\mu_{0}|\|u^{n}\|_{B^{s}_{p,r}}+|\mu_{0}|\|u^{n}-u\|_{B^{s-1}_{p,r}}\\
\leq \ &
2M|\mu_{0}^{n}-\mu_{0}|+|\mu_{0}|\|u^{n}-u\|_{B^{s-1}_{p,r}},
\end{align*}
$$\|\partial_{x}P(D)(\frac{1}{2}(u^{n}_{x})^{2}-\frac{1}{2}u_{x}^{2})\|_{B^{s-1}_{p,r}}\leq
\frac{1}{2}\|(u^{n}_{x})^{2}-u_{x}^{2}\|_{B^{s-1}_{p,r}}\leq M
\|u^{n}_{x}-u_{x}\|_{B^{s-1}_{p,r}},$$
\begin{align*}
\|\partial_{x}P(D)(\frac{1}{2}(\rho^{n})^{2}-\frac{1}{2}\rho^{2})\|_{B^{s-1}_{p,r}}
\leq \ &\frac{1}{2}
\|(\rho^{n})^{2}-\rho^{2}\|_{B^{s-1}_{p,r}}\\
\leq \
&\frac{1}{2}\|\rho^{n}+\rho\|_{B^{s-1}_{p,r}}\|\rho^{n}-\rho\|_{B^{s-1}_{p,r}}\\
\leq \ & M\|\rho^{n}-\rho\|_{B^{s-1}_{p,r}}.
\end{align*}
It then follows that
\begin{align}
&\nonumber\|F^{n}-F\|_{B^{s-1}_{p,r}}\\
\nonumber \leq & \
3M\|u^{n}_{x}-u_{x}\|_{B^{s-1}_{p,r}}+M\|\rho^{n}-\rho\|_{B^{s-1}_{p,r}}\\
\nonumber & \ \ \ \ \ \ \ \  \ \ \ \ \  \ \ \ \ \
+2M|\mu_{0}^{n}-\mu_{0}|+|\mu_{0}|\|u^{n}-u\|_{B^{s-1}_{p,r}}\\
\nonumber \leq & \
3M\|v^{n}\|_{B^{s-1}_{p,r}}+3M\|w^{n}-u_{x}\|_{B^{s-1}_{p,r}}+M\|g^{n}\|_{B^{s-1}_{p,r}}\\
 & \ \ \ \ \  \ \
+M\|f^{n}-\rho\|_{B^{s-1}_{p,r}}+2M|\mu_{0}^{n}-\mu_{0}|+|\mu_{0}|\|u^{n}-u\|_{B^{s-1}_{p,r}}.
\end{align}
Moreover,
\begin{align}
&\nonumber\|u^{n}_{x}\rho^{n}-u_{x}\rho\|_{B^{s-1}_{p,r}}\\
\nonumber\leq & \
\|(u^{n}_{x}-u_{x})\rho^{n}\|_{B^{s-1}_{p,r}}+\|u_{x}(\rho^{n}-\rho)\|_{B^{s-1}_{p,r}}\\
\nonumber \leq & \
M\|u^{n}_{x}-u_{x}\|_{B^{s-1}_{p,r}}+M\|\rho^{n}-\rho\|_{B^{s-1}_{p,r}}\\
 \leq & \
M\|v^{n}\|_{B^{s-1}_{p,r}}+M\|w^{n}-u_{x}\|_{B^{s-1}_{p,r}}+M\|g^{n}\|_{B^{s-1}_{p,r}}+M\|f^{n}-\rho\|_{B^{s-1}_{p,r}}.
\end{align}
Combining (3.21)-(3.24), we get
\begin{align*}
&\|v^{n}\|_{B^{s-1}_{p,r}}+\|g^{n}\|_{B^{s-1}_{p,r}}\leq
e^{CMT}\{\|u^{n}_{0x}-u_{0x}\|_{B^{s-1}_{p,r}}+\|\rho_{0}^{n}-\rho_{0}\|_{B^{s-1}_{p,r}}\\
&\ \ \ \ \  \ \ \ \ \ \  \ \ \ \ \ \ \ \ \ +\int_{0}^{t}[4M(\|v^{n}\|_{B^{s-1}_{p,r}}+\|g^{n}\|_{B^{s-1}_{p,r}})+4M\|w^{n}-u_{x}\|_{B^{s-1}_{p,r}}\\
&\ \ \ \ \ \ \ \ \ \  \ \ \ \ \ \ \ \ \ \
+2M\|f^{n}-\rho\|_{B^{s-1}_{p,r}}+2M|\mu_{0}^{n}-\mu_{0}|+|\mu_{0}|\|u^{n}-u\|_{B^{s-1}_{p,r}}]d\tau\}.
\end{align*}
Note that $$z_{0}^{n}\rightarrow z_{0}^{\infty}:=z_{0} \ \ \text{in}
\ \ B^{s}_{p,r}\times B^{s-1}_{p,r},$$
$$w^{n}\rightarrow w=u_{x}
\ \ \text{in} \ \ C([0,T]; B^{s-1}_{p,r})\ \ \ , \ \ \
f^{n}\rightarrow \rho \ \ \text{in} \ \ C([0,T]; B^{s-1}_{p,r}),$$
and $$u^{n}\rightarrow u \ \ \text{in} \ C([0,T]; B^{s-1}_{p,r}).$$
By the Lebesgue dominated convergence theorem, we have
$$\lim\limits_{n\rightarrow
\infty}(\|v^{n}\|_{B^{s-1}_{p,r}}+\|g^{n}\|_{B^{s-1}_{p,r}})\leq
4Me^{CMT}\int_{0}^{t}\lim\limits_{n\rightarrow
\infty}(\|v^{n}\|_{B^{s-1}_{p,r}}+\|g^{n}\|_{B^{s-1}_{p,r}})d\tau.$$
Applying Gronwall's inequality, we have
$$v^{n}\rightarrow 0 \ \ \text{in} \ \
C([0,T]; B^{s-1}_{p,r})\ \ \ \text{and} \ \ \ g^{n}\rightarrow 0 \ \
\text{in} \ \ C([0,T]; B^{s-1}_{p,r}),$$ which implies
$$u^{n}_{x}\rightarrow u_{x} \ \ \text{in} \ \ C([0,T];
B^{s-1}_{p,r})\ \ \ \text{and} \ \ \ \rho^{n}\rightarrow \rho \ \
\text{in} \ \ C([0,T]; B^{s-1}_{p,r}).$$ This completes the proof of
Theorem 3.1.

Since the Sobolev space $H^{s}=B_{2,2}^{s},$ Theorem 3.1 implies
that if $(u_{0}, \rho_{0})\in H^{s}\times H^{s-1}$ with
$s>\frac{3}{2}$ and $s\neq \frac{5}{2},$ we can obtain the local
well-posedness to the system (3.2) in $H^{s}\times H^{s-1}$ with
$\frac{3}{2}<s \neq \frac{5}{2}$. Combining the corresponding local
well-posedness result in \cite{l-y1} (where $s\geq 2$ is obtained)
and letting $p=r=2$ in Theorem 3.1, we get the following main result
of this section:\\
\newline
\textbf{Theorem 3.2.} Given $z_{0}=(u_{0},\rho_{0})\in H^{s}\times
H^{s-1}, s> \frac{3}{2},$ there exists a maximal $T=T(\parallel
z_{0}\parallel_{H^{s}\times H^{s-1}})>0$, and a unique solution
$z=(u,\rho)$ to the system (3.2)(or (1.2)) such that
$$
z=z(\cdot,z_{0})\in C([0,T); H^{s}\times H^{s-1})\cap
C^{1}([0,T);H^{s-1}\times H^{s-2}).
$$
Moreover, the solution depends continuously on the initial data,
i.e. the mapping
$$z_{0}\rightarrow z(\cdot,z_{0}): H^{s}\times H^{s-1}\rightarrow
C([0,T); H^{s}\times H^{s-1})\cap C^{1}([0,T);H^{s-1}\times H^{s-2})
$$
is continuous.

\section{The precise blow-up scenario}

In this section, we present the precise blow-up scenario and global
existence for solutions to the system (3.2) in Sobolev spaces.\\
\newline
\textbf{Lemma 4.1} (\cite{ fu, constantin}).\ If $f\in
H^{1}(\mathbb{S})$ is such that $\int_{\mathbb{S}}f(x)dx=0,$ then we
have
$$\max\limits_{x\in\mathbb{S}}f^{2}(x)\leq
\frac{1}{12}\int_{\mathbb{S}}f_{x}^{2}(x)dx.$$

Note that $\int_{\mathbb{S}}(u(t,x)-\mu_{0})dx=\mu_{0}-\mu_{0}=0.$
By Lemma 4.1, we find that
$$\max\limits_{x\in\mathbb{S}}[u(t,x)-\mu_{0}]^{2}\leq
 \frac{1}{12}\int_{\mathbb{S}}u_{x}^{2}(t,x)dx\leq
 \frac{1}{12}\mu_{1}^{2}.$$ So we have
 \begin{equation}
 \|u(t,\cdot)\|_{L^{\infty}(\mathbb{S})}\leq
 |\mu_{0}|+\frac{\sqrt{3}}{6}\mu_{1}.
 \end{equation}

Consider now the following initial value problem

\begin{equation}
\left\{\begin{array}{ll}q_{t}=u(t,-q)+2\gamma_{2},\ \ \ \ t\in[0,T), \\
q(0,x)=x,\ \ \ \ x\in\mathbb{R}, \end{array}\right.
\end{equation}
where $u$ denotes the first component of the solution $z$ to the
system (3.2). Then we have the following two useful lemmas.

Similar to the proof of Lemma 4.1 in \cite{Y1}, applying classical
results in the theory of ordinary differential equations, one can
obtain the following result on $q$ which is crucial in the proof of
blow-up scenarios.\\
\newline
\textbf{Lemma 4.2} (\cite{l-y1}).\ Let $u\in C([0,T); H^{s})\bigcap
C^{1}([0,T); H^{s-1}), s >\frac{3}{2}$. Then Eq.(4.2) has a unique
solution $q\in C^1([0,T)\times \mathbb{R};\mathbb{R})$. Moreover,
the map $q(t,\cdot)$ is an increasing diffeomorphism of $\mathbb{R}$
with
$$
q_{x}(t,x)=exp\left(-\int_{0}^{t}u_{x}(s,-q(s,x))ds\right)>0, \ \
(t,x)\in [0,T)\times \mathbb{R}.$$\\
\newline
\textbf{Lemma 4.3} (\cite{l-y1}). \ Let
$z_{0}=\left(\begin{array}{c}
                                u_{0} \\
                                \rho_{0} \\
                              \end{array}
                            \right)
\in H^{s}\times H^{s-1}$, $s >\frac{3}{2}$ and let $T>0$ be the
maximal existence time of the corresponding solution
                  $z=\left(\begin{array}{c}
                                      u \\
                                      \rho \\
                                    \end{array}
                                  \right)$
to (3.2). Then we have
\begin{equation}
\rho(t,-q(t,x))q_{x}(t,x)=\rho_{0}(-x), \ \ \ \forall \
(t,x)\in[0,T)\times \mathbb{S}.
\end{equation}
Moreover, if there exists $M> 0$ such that $u_{x}\leq M$ for all
$(t,x)\in [0,T)\times \mathbb{S}$, then
$$\|\rho(t,\cdot)\|_{L^{\infty}}\leq
e^{MT}\|\rho_{0}(\cdot)\|_{L^{\infty}}, \ \ \ \forall \ t\in[0,
T).$$

Our next result implies that the wave breaking to the system (3.2)
is determined only by the slope of $u$ but not the slope of $\rho.$\\
\newline
\textbf{Theorem 4.1.} Let $z_{0}=\left(
                                                     \begin{array}{c}
                                                       u_{0} \\
                                                       \rho_{0} \\
                                                     \end{array}
                                                   \right)
\in H^s\times H^{s-1}$ with $s>\frac 3 2$ and let $T$ be the maximal
existence time of the solution $z=\left(
                                                     \begin{array}{c}
                                                       u \\
                                                      \rho \\
                                                     \end{array}
                                                   \right) $ to the system (3.2), which is guaranteed by Theorem 3.2.
If $T<\infty$, then $$\int_0^T \|\partial_x u(\tau)\|_{L^\infty}
d\tau=\infty.$$ \

The proof of the theorem is similar to that of Theorem 4.1 in
\cite{gui2,yan}, so we omit it here.\

Next, we first recall a useful lemma before giving the next
result.\\
\newline
\textbf{Lemma 4.4} (\cite{Constantin 4}).\ Let $t_{0}>0$ and $v\in
C^{1}([0,t_{0}); H^{2}(\mathbb{R}))$. Then for every $t\in[0,t_{0})$
there exists at least one point $\xi(t)\in \mathbb{R}$ with
$$ m(t):=\inf_{x\in \mathbb{R}}\{v_{x}(t,x)\}=v_{x}(t,\xi(t)),$$ and
the function $m$ is almost everywhere differentiable on $(0,t_{0})$
with $$ \frac{d}{dt}m(t)=v_{tx}(t,\xi(t)) \ \ \ \ a.e.\ on \
(0,t_{0}).$$
\newline
\textbf{Remark 4.1.} If $v\in C^{1}([0,t_{0}); H^{s}(\mathbb{R})), \
s>\frac{3}{2},$ then Lemma 4.4 also holds true. Meanwhile, Lemma 4.4
works analogously for
$$M(t):= \sup\limits_{x\in\mathbb{R}}\{{v_{x}(t,x)}\}.$$

Our next result describes the precise blow-up scenario for
sufficiently regular solutions to the system (3.2). This result for
the system improve considerably the earlier result in \cite{l-y1}.\\
\newline
\textbf{Theorem 4.2.} Let $z_{0}=\left(
                                                     \begin{array}{c}
                                                       u_{0} \\
                                                       \rho_{0} \\
                                                     \end{array}
                                                   \right)
\in H^s\times H^{s-1}, s>\frac{5}{2},$ and let T be the maximal
existence time of the solution $z=\left(
                                    \begin{array}{c}
                                      u \\
                                      \rho \\
                                    \end{array}
                                  \right)
$ to the system (3.2) with the initial $z_{0}$. Then the the maximal
existence time $T$ is finite if and only if
$$\limsup\limits_{t\rightarrow
T}\sup\limits_{x\in\mathbb{S}}\{u_{x}(t,x)\}=+\infty.$$
\newline
\textbf{Proof}
On one hand, by Sobolev's imbedding theorem it is clear that if
$$\limsup\limits_{t\rightarrow
T}\sup\limits_{x\in\mathbb{S}}\{u_{x}(t,x)\}=+\infty,$$ then the the
maximal existence time $T<\infty$.

On the other hand, assume that the maximal existence time $T$ is
finite and there exists a $M>0$ such that
\begin{eqnarray}
u_{x}(t,x)\leq M,\ \ \ \forall(t,x)\in [0,T)\times\mathbb{S}.
\end{eqnarray}
Then, from Lemma 4.3, we have
\begin{equation}
\|\rho(t,\cdot)\|_{L^{\infty}}\leq e^{MT}\|\rho_{0}\|_{L^{\infty}},
\ \ \ \forall \ t\in[0, T).
\end{equation}
Let $m(t)=\min\limits_{x\in \mathbb{S}}\{u_{x}(t,x)\}.$ It follows
from Remark 4.1 that there is a point $(t,\xi(t))\in [0,T)\times
\mathbb{S}$ such that $m(t)=u_{x}(t,\xi(t)).$ Moreover,
$u_{xx}(t,\xi(t))=0.$ Evaluating (3.3) on $(t, \xi(t))$ we get
\begin{align*}
\frac{d}{dt}m(t)&=-2\mu_{0}u(t,
\xi(t))+\frac{1}{2}m^{2}(t)-\frac{1}{2}\rho^{2}(t,
\xi(t))+a\\
&\geq-2\mu_{0}u(t, \xi(t))-\frac{1}{2}\rho^{2}(t, \xi(t)).
\end{align*}
By (4.1) and (4.5), we get $$\frac{d}{dt}m(t)\geq
-2|\mu_{0}|(|\mu_{0}|+\frac{\sqrt{3}}{6}\mu_{1})-\frac{1}{2}\left(e^{MT}\|\rho_{0}\|_{L^{\infty}}\right)^{2}.$$
Integrating this inequality on $(0,t)$, we have
$$m(t)\geq m(0)-\left(2|\mu_{0}|(|\mu_{0}|+\frac{\sqrt{3}}{6}\mu_{1})
+\frac{1}{2}\left(e^{MT}\|\rho_{0}\|_{L^{\infty}}\right)^{2}\right)T.$$
That is, $$\min\limits_{x\in \mathbb{S}}u_{x}(t,x)\geq
\min\limits_{x\in
\mathbb{S}}u_{0}^{\prime}(x)-\left(2|\mu_{0}|(|\mu_{0}|+\frac{\sqrt{3}}{6}\mu_{1})
+\frac{1}{2}\left(e^{MT}\|\rho_{0}\|_{L^{\infty}}\right)^{2}\right)T,$$
which together with (4.4) and $T<\infty$ implies that
$$\int_0^T \|\partial_x u(\tau)\|_{L^\infty} d\tau<\infty.$$
This contradicts Theorem 4.1.

Furthermore, if $\gamma_{1}=2\gamma_{2},$ then we get the following
sharper conclusion for $s.$\\
\newline
\textbf{Theorem 4.3.} Let $z_{0}=\left(
                                                     \begin{array}{c}
                                                       u_{0} \\
                                                       \rho_{0} \\
                                                     \end{array}
                                                   \right)
\in H^s\times H^{s-1}, s>\frac{3}{2},$ and let T be the maximal
existence time of the solution $z=\left(
                                    \begin{array}{c}
                                      u \\
                                      \rho \\
                                    \end{array}
                                  \right)
$ to the system (3.2) with the initial $z_{0}$. Assume
$\gamma_{1}=2\gamma_{2},$ then the the maximal existence time $T$ is
finite if and only if
$$\limsup\limits_{t\rightarrow
T}\sup\limits_{x\in\mathbb{S}}\{u_{x}(t,x)\}=+\infty.$$
\newline
\textbf{Proof}
On one hand, by Sobolev's imbedding theorem it is clear that if
$$\limsup\limits_{t\rightarrow
T}\sup\limits_{x\in\mathbb{S}}\{u_{x}(t,x)\}=+\infty,$$ then the the
maximal existence time $T<\infty$.

On the other hand, assume that the the maximal existence time $T$ is
finite and there exists a $M>0$ such that
\begin{eqnarray}
u_{x}(t,x)\leq M,\ \ \ \forall(t,x)\in [0,T)\times\mathbb{S}.
\end{eqnarray}
Then, from Lemma 4.3, we have
\begin{equation}
\|\rho(t,\cdot)\|_{L^{\infty}}\leq e^{MT}\|\rho_{0}\|_{L^{\infty}},
\ \ \ \forall \ t\in[0, T).
\end{equation}
By (4.2) and the condition $\gamma_{1}=2\gamma_{2},$ we have
\begin{align}
\nonumber\frac{du_x(t,-q(t,x))}{dt}
&=u_{xt}(t,-q(t,x))-u_{xx}(t,q(t,x))q_t(t,x)\\
&=(u_{tx}-(u+\gamma_{1})u_{xx})(t,-q(t,x)).
\end{align}
Evaluating (3.3) on $(t, -q(t,x))$ we get
\begin{align*}
&\frac{du_x(t,-q(t,x))}{dt}\\
=&-2\mu_{0}u(t,
-q(t,x))+\frac{1}{2}u_{x}^{2}(t,-q(t,x))-\frac{1}{2}\rho^{2}(t,
-q(t,x))+a\\
\geq&-2\mu_{0}u(t, -q(t,x))-\frac{1}{2}\rho^{2}(t,-q(t,x)).
\end{align*}
Similar to the proof of Theorem 4.2, we obtain
$$\min\limits_{x\in \mathbb{S}}u_{x}(t,x)\geq
\min\limits_{x\in
\mathbb{S}}u_{0}^{\prime}(x)-\left(2|\mu_{0}|(|\mu_{0}|+\frac{\sqrt{3}}{6}\mu_{1})
+\frac{1}{2}\left(e^{MT}\|\rho_{0}\|_{L^{\infty}}\right)^{2}\right)T,$$
which together with (4.6) and $T<\infty$ implies that
$$\int_0^T \|\partial_x u(\tau)\|_{L^\infty} d\tau<\infty.$$
This contradicts Theorem 4.1.

Next, we state an improved global existence theorem, which improves
the result of the global solutions in \cite{l-y1}, where the special
case $s=2$ is required. However, the proof of this improved result
is the same as the proof of corresponding result in \cite{l-y1}, so
we omit it here.\\
\newline
\textbf{Theorem 4.4.} Let $z_{0}=\left(
                                                     \begin{array}{c}
                                                       u_{0} \\
                                                       \rho_{0} \\
                                                     \end{array}
                                                   \right)
\in H^s\times H^{s-1}, s>\frac{3}{2} ,$  and T be the maximal time
of the solution $z=\left(
                                    \begin{array}{c}
                                      u \\
                                      \rho \\
                                    \end{array}
                                  \right)
$ to the system (3.2) with the initial data $z_0$. If
$\gamma_{1}=2\gamma_{2}$ and $\rho_{0}(x)\neq 0$ for all
$x\in\mathbb{S}$, then the corresponding solution $z$ exists
globally in time.

\end{document}